\documentclass[final,3p,12pt]{elsarticle} 

\usepackage{graphicx}
\usepackage{subfigure}
\usepackage{epstopdf}
\usepackage{siunitx}
\usepackage{gensymb}

\usepackage{amssymb}

\usepackage{amsmath}

\usepackage{lineno}
\usepackage{marginnote}
\usepackage{color,soul}
\definecolor{G}{rgb}{.80,0.0,.0}
\definecolor{A}{rgb}{.0,.40,.90}
\usepackage{setspace}
\usepackage{multirow}
\usepackage{url}
\usepackage{ctable}
\usepackage{bm}	
\usepackage{enumerate}

\journal{{Computer Methods in Applied Mechanics and Engineering}}

\usepackage{siunitx}
\usepackage{textcomp}

\usepackage{xcolor}

\begin{document}

\begin{frontmatter}


\title{
An enhanced VEM formulation for plane elasticity
}

\author{A.M. D'Altri$^{1}$}
\author{S. de Miranda$^{1}$\corref{cor1}}
\author{L. Patruno$^{1}$}
\author{E. Sacco$^{2}$} 
\cortext[cor1]{corresponding author: stefano.demiranda@unibo.it}


\address{%

$^{1}$ Department of Civil, Chemical, Environmental, and Materials Engineering (DICAM), University of Bologna, Viale del Risorgimento 2, 40136 Bologna, Italy

$^{2}$ Department of Structures for Engineering and Architecture, University of Naples Federico II, Naples, Italy

}




\singlespacing
\begin{abstract} 
\small
  In this paper, an enhanced Virtual Element Method (VEM) formulation is proposed for plane elasticity.
It is based on the improvement of the strain representation within the element, without altering the degree of the displacement interpolating functions on the element boundary. 
The idea is to fully exploit polygonal elements with a high number of sides, a peculiar VEM feature, characterized by 
many displacement degrees of freedom on \color{black}
 the element boundary, even if a low interpolation order is assumed over each side.
The proposed approach is framed within a generalization of the classic VEM formulation, obtained by introducing an energy norm in the projection operator definition.
Although such generalization may mainly appear to have a formal value,
it allows to effectively point out the mechanical meaning of the  quantities involved in the projection operator definition and to drive the selection of the enhanced representations.
Various enhancements are proposed and tested through several numerical examples. 
Numerical results successfully show the capability of the enhanced VEM formulation to (i) considerably increase accuracy (with respect to standard VEM) while keeping the optimal convergence rate, (ii) bypass the need of stabilization terms in many practical cases, (iii) obtain natural serendipity elements  in many practical cases, and (vi) 
effectively treat also nearly incompressible materials.
\end{abstract}


\end{frontmatter}

\onehalfspacing

\section{Introduction}

The Virtual Element Method (VEM) is a recent and interesting numerical technique for solving boundary value problems governed by a system of partial differential equations. 
It has been proposed by Beir\~ao da Veiga et al. \cite{beirao2013basic}, who founded the method and called it  ``Virtual Element Method''. 
The idea is to modify (improve) the construction of the standard finite element method to obtain a method that is at same time simple and able to preserve the polynomial accuracy, accounting for the specific features of the mimetic finite difference method \cite{hyman1999orthogonal,beirao2010mimetic}  and of the polygonal finite element method \cite{sukumar2004}.
The method has been, then, developed by the VEM proposing research team to solve different types of partial differential equations governing interesting problems in Mechanics, such as  linear elasticity \cite{beirao2013virtual},   Kirchhoff-Love plate bending problems \cite{brezzi2013},  parabolic problems considering time‐dependent diffusion equations.
A sort of manual (in mathematical sense) of the method is available in \cite{beirao2014hitchhiker}.

In the computational mechanics community, the method arouse great interest.
In fact, several studies concerning the accuracy, the computational effectiveness and the applicability of the method to different mechanical problems have been developed in the last years, proposing even some modifications and improvements \cite{beirao2016serendipity,wriggers2017low,Artioli_2017,artioli2018family,dassi2019bricks,artioli2020curvilinear,dassi2020three,mascotto2019nonconforming,mora2020virtual,de2019serendipity,artioli2019dual,chi2017some,park2020numerical}.
The advantages of the method lie in the fact that it is easy to implement, does not require as in Finite Element Method (FEM) the definition of the parent element, the approximation order can be quite simply modified.
Of course, the greatest advantage is that in two dimensional problems the element can have how many sides as required in order to construct a suitable mesh for the problem to solve.
All these advantages have been exploited for example in contact problems \cite{wriggers2016virtual}, and in fracture mechanics \cite{nguyen2018virtual,aldakheel2018phase,hussein2019computational,artioli2020vem}.

In solid mechanics, the crucial point of the construction of the method consists in the definition of the approximation of the displacement field only on the boundary of the element, without an explicit representation of the displacements inside the element. For this reason the strain is computed by a projection operation and a stabilization term is required to avoid the insurgence of zero-energy modes. 

{According to the standard displacement-based VEM formulation \cite{artioli2017arbitrary}, this projection operation is conducted assuming
	 the strain polynomial approximation one degree lower than the one used for displacements interpolation on the element boundary.
This assumption is made independently from the number of sides of the element \cite{artioli2017arbitrary}.
Therefore, this assumption can lead to paradoxical situations in elements with high number of sides, given that
the many degrees of freedom
available on the element boundary are not fully exploited to compute the strain.

An example of this paradoxical situation is reported in
 \cite{artioli2018high}, where} the micromechanical analysis of a unit cell of a periodic composite material has been performed adopting the VEM technique. 
Accounting that an inclusion in a infinite medium subjected to prescribed average strain is characterized by uniform strain, the fiber of the repetitive unit cell could have an almost uniform strain distribution, depending on the inclusion volume matrix. 
For this reason, the core of the inclusion, at a certain distance from the fiber-matrix interface, has been discretized with only one virtual element, saving computational efforts. This element is characterized by 140 sides and 140 nodes, as linear interpolation functions have been chosen for the displacements at the boundary. 
Thus, constant strain is considered inside the element. Of course, to get constant strain much less nodes are required, so that the element is mostly stabilized and all the nodal unknowns are not much used to significantly improve the strain field inside the element and hence the solution of the problem.
Moreover, even the stabilization techniques are not always completely reliable, as discussed in \cite{wriggers2017efficient}.

Indeed, it is clear that even leaving a given degree of the polynomial approximating the displacement field on the boundary, the number of nodes (and hence sides) of a virtual element can lead to an improvement of the solution. This point is discussed  in \cite{artioli2019equilibrium,daltri2020}, where the recovery of the stress field inside the element is evaluated as post processing of the solution for different number of sides.

Moreover, a very recent contribution  in this field has been given in \cite{LOURENCO2020sharper}, where a sharper error analysis has been developed for the VEM to separate the element boundary and the element interior contributions to the error, in order to yield  more accurate discrete solutions.

\color{black}

In this paper, an enhanced VEM formulation is proposed for plane elasticity.
Basically, the proposed enhancement is based on the improvement of the strain representation within the element, without altering the degree of the displacement interpolating functions on the element boundary. 
Accordingly, this approach allows to fully exploit elements with a high number of sides, a peculiar VEM feature, characterized by 
 many displacement degrees of freedom on \color{black}
the element boundary 
(i.e. a lot of information resides on the element boundary in its entirety), 
even if a low interpolation degree is assumed over each side.

The proposed approach originates within the framework of a generalization of the classic VEM formulation \cite{artioli2017arbitrary}, obtained by introducing an energy norm in the projection operator definition.
Although such generalization may mainly appear to have a formal value, 
it allows to effectively point out the mechanical meaning of the  quantities involved in the projection operator definition and to drive the selection of the enhanced representations.

%

Various enhancements are proposed 
and tested through several numerical examples to check their performance with respect to standard displacement-based VEM. 
Particularly, the enhanced VEM performance in terms of convergence curves, the possibility to bypass the need of stabilization terms, and the performance  in nearly incompressible materials are analyzed and discussed.

The paper is organized as follows. Section 2 summarizes the basic equations of the plane elasticity problem. Section 3 recalls the standard VEM formulation. Section 4 introduces and describes the enhanced VEM formulation. Section 5 collects and discusses the numerical results. Finally, Section 6 highlights the conclusions of this research.

\begin{table}[htb]

	\footnotesize
	\centering
	\begin{tabular*}{\textwidth}{|c @{\extracolsep{\fill}} p{6cm} c p{6cm}|}
		\cline{1-4} 
		\multicolumn{4}{|l|}{ }    \\
		\multicolumn{4}{|l|}{\textbf{List of main symbols} }   \\
		\multicolumn{4}{|l|}{ }    \\
		
		$\mathbf{b}$  & distributed volume forces&$s$  & stabilization order ($s=k$ in standard VEM)\\
		$\mathbf{C}$   & elasticity matrix &$\mathbf{U}$ & space of admissible displacements\\
		$k$ &  polynomial degree of approximation of the displacement on the element boundary &$\mathbf{U}_0$   & space of the variations of admissible displacements \\
		$\mathbf{K}$ & total local stiffness matrix &$\mathbf{v}$ &  displacement vector  \\
		$\mathbf{K}_c$ & consistent part of the stiffness matrix &$\mathbf{v}_\mathit{h}$  & virtual displacement vector  \\
		$\mathbf{K}_s$  &  stabilization part of the stiffness matrix&$\tilde{\mathbf{v}}_\mathit{h}$ & vector of approximated displacement  on the element boundary\\
		$\mathbf{L}$   & compatibility operator &		$\tilde{\mathbf{V}}$ &  vector collecting the displacements of the nodes on the element boundary \\
		$\mathbf{L}^\text{T}$  & equilibrium operator&		$\hat{\mathbf{V}}$  & vector collecting the  moments of the virtual displacement\\
		$m$  & number of vertexes of the element&		$\bm{\varepsilon}$  & strain vector\\
		$n$  &  total number of the element degrees of freedom ($ n = 2mk + p (p + 1) $ for standard VEM)&		${\bm{{\varepsilon}}}^P$ & projected strain vector \\
		$\mathbf{N}_E$  & matrix containing the outward unit normal vector direction cosines on $\partial \mathit{\Omega}_E$  &		$\bm{\sigma}$  & stress vector \\		
		$\mathbf{N}^P$  & matrix of the approximation functions for ${\bm{{\varepsilon}}}^P$&		$\mathit{\Omega}$ & domain of the  body\\
		$\breve{\mathbf{N}}^P$  &  matrix of the approximation functions for $\mathbf{C}{\bm{{\varepsilon}}}^P$ &		$\mathit{\Omega_E}$ & domain of the typical  element \\
		$\mathbf{N}^V$  & matrix of the  approximation functions of $\tilde{\mathbf{v}}_\mathit{h}$ on the element boundary&		$\partial \mathit{\Omega}$  & boundary of the body \\
		$p$ &  polynomial degree of the approximation functions in $\mathbf{N}^P$ or $\breve{\mathbf{N}}^P$   ($p=k-1$ in standard VEM) &		$\partial \mathit{\Omega_E}$  & boundary of the typical element   \\
		
		\cline{1-4}	
	\end{tabular*}
\end{table}

\color{black}

\section{Basic equations}\label{problem}
Consider a body that occupies a closed and bounded region in $\mathbb{R}^2$, whose boundary and open set are denoted by $\partial \mathit{\Omega}$ and $\mathit{\Omega}$, respectively. 
The generic configuration of the body is described by the displacement vector $\mathbf{v}$, and by the
strain and stress tensors, whose components are arranged in form of 3-components vectors
$\bm{\varepsilon}$ and $\bm{\sigma}$, respectively. The distributed volume forces $\mathbf{b}$ are prescibed in $\mathit{\Omega}$.

The governing equations of the linear elastic problem in $\mathit{\Omega}$ (i.e. the strain-displacement relationship, the equilibrium equation and the constitutive law) can be written as:
\begin{equation}
\bm{\varepsilon}=\mathbf{Lv},
\end{equation}
\begin{equation}
\mathbf{L}^\text{T}\bm{\sigma}+\mathbf{b}=\mathbf{0},
\end{equation}
\begin{equation}
\bm{\sigma}=\mathbf{C} \bm{\varepsilon},
\end{equation}
where $\mathbf{L}$ and $\mathbf{L}^\text{T}$ are the compatibility and equilibrium operators and $\mathbf{C}$ the elasticity matrix.
The above governing equations are completed by the boundary conditions on $\partial \mathit{\Omega}$.


The variational formulation of the linear elastic problem in terms of displacements stems from the virtual work principle:
\begin{equation} \label{PLV}
\begin{cases}
\mathrm{Find} \: \mathbf{v}\in\mathbf{U} \: \mathrm{such \: that} \\
{\displaystyle \int_{\mathit{\Omega}} (\mathbf{L} \, \delta \mathbf{v})^\text{T}\mathbf{C}\mathbf{Lv}  \; \text{dA}}
=
{\displaystyle \int_{\mathit{\Omega}} \delta \mathbf{v}^\text{T}\mathbf{b} \; \text{dA} }
\;\;\;\; \forall \: \delta \mathbf{v} \in \mathbf{U}_0
\end{cases}
\end{equation}
where $\mathbf{U}$ is the space of the admissible displacements and $\mathbf{U}_0$ is the space of the variations of admissible displacements.

\section{Standard VEM formulation}\label{formulation}
In this section, the standard displacement-based VEM formulation is briefly summarized. The interested reader can refer to \cite{artioli2017arbitrary} for the details of the implementation of the standard displacement-based VEM in a computer code. 

To solve  the problem in Eq. (\ref{PLV}) through the VEM scheme, the domain $\mathit{\Omega}$ is discretized by means of nonoverlapping polygons with straight edges. Each one of such polygons can be denoted with $\mathit{\Omega_E}$, while $\partial \mathit{\Omega_E}$ indicates its boundary. 
In the VEM, the approximated displacement field in the element interior $\mathit{\Omega_E}$, denoted as  $\mathbf{v}_\mathit{h}$ in the following, is assumed to be not explicitly known and, thus, is also referred to as virtual. On the contrary, the approximated displacement field on the element boundary $\partial \mathit{\Omega_E}$, denoted in the following by 
$ \tilde{\mathbf{v}}_\mathit{h}$, is assumed to be explicitly known and written as:
\begin{equation} \label{eq_virtual-displacement}
\tilde{\mathbf{v}}_\mathit{h} = \mathbf{N}^V \tilde{\mathbf{V}},
\end{equation}  
where $\mathbf{N}^V$ is the matrix of the  approximation functions on $\partial \mathit{\Omega_E}$  and $\tilde{\mathbf{V}}$ is the vector collecting the displacements of the nodes on the element boundary. 
Approximated displacements $ \mathbf{v}_\mathit{h}$ and $ \tilde{\mathbf{v}}_\mathit{h}$ coincide on $\partial \mathit{\Omega_E}$ (i.e. $\tilde{\mathbf{v}}_\mathit{h}$ is the restriction of $\mathbf{v}_\mathit{h}$ to $\partial \mathit{\Omega_E}$).
Vector $\tilde{\mathbf{V}}$ has $2mk$ components, where $m$ is the number of vertexes of the polygonal element and $k$ the degree of the polynomial representation assumed over the element boundary.




\subsection{Projection operator and consistent term}\label{projection}
As neither $\mathbf{v}_\mathit{h}$ nor its gradient are explicitly
computable in the element interior, the key aspect of the VEM consists in the introduction of a projection operator, which represents the approximated strain associated with the virtual displacement,
called projected strain and denoted by  ${\bm{{\varepsilon}}}^P$ in the following.
Accordingly, given the virtual displacement $\mathbf{v}_\mathit{h}$,
 $\bm{{\varepsilon}}^P$ can be defined as the unique function
 that minimizes
\begin{equation} \label{eq_requirement-projector}
\left\| \bm{{\varepsilon}}^P - \bm{{\varepsilon}}(\mathbf{v}_\mathit{h})\right\|_{norm}
\end{equation}
where subscript $norm$ indicates the type of norm used. 
In order to perform the miminization of Eq. (\ref{eq_requirement-projector}), the following assumptions have to be made:
\begin{enumerate}[(i)]
	\item the norm to be used;
	\item the  representation to be used for approximated quantities.
\end{enumerate}

As regards the first point, in the standard VEM, an  $\mathit{L}^2$ norm is used 
and, thus, the minimization of  Eq.  (\ref{eq_requirement-projector})  leads to: find $\bm{{\varepsilon}}^P \in \mathsf{P}_p(\Omega_E)$ such that
\begin{equation} \label{min_L2_a}
 \int_{\mathit{\Omega}_E} [\bm{{\varepsilon}}^P-\bm{\varepsilon}(\mathbf{v}_\mathit{h})]^\text{T} \delta {\bm{{\varepsilon}}}^P \text{dA} = 0\; \; \; \forall \;\delta {\bm{{\varepsilon}}}^P\in \mathsf{P}_p(\Omega_E).
\end{equation}

\color{black}
As regards the second point, in the standard VEM, an uncoupled polynomial approximation of each strain component complete up to degree $p$ is assumed  for $\bm{{\varepsilon}}^P$,
with the degree $p$ linked to the order  $k$  used in $\mathbf{N}^V$ by the relationship:
\begin{equation} \label{p_and_k}
p=k-1,
\end{equation}
\noindent that is:
\begin{equation} \label{eq_projector-operator}
\bm{{\varepsilon}}^P= \mathbf{N}^P\hat{\bm{{\varepsilon}}},
\end{equation}
\setcounter{MaxMatrixCols}{30}
\begin{equation}  \label{eq_NP_div1}
\mathbf{N}^P = 
\begin{bmatrix}
1 & 0 & 0 & x & 0 & 0 & y & 0 & 0 & \cdots & y^p & 0 & 0 \\
0 & 1 & 0 & 0 & x & 0 & 0 & y & 0 & \cdots & 0 & y^p & 0\\
0 & 0 & 1 & 0 & 0 & x & 0 & 0 & y & \cdots & 0& 0 & y^p
\end{bmatrix}.
\end{equation}



\noindent Using Eq. (\ref{eq_projector-operator}), Eq. (\ref{min_L2_a}) can be rewritten as:
\begin{equation} \label{min_L2}
\delta \hat{\bm{{\varepsilon}}}^\text{T} \int_{\mathit{\Omega}_E} ( \mathbf{N}^P)^\text{T} [ \mathbf{N}^P\hat{\bm{{\varepsilon}}} - \bm{{\varepsilon}}(\mathbf{v}_\mathit{h})] \text{dA} = 0\; \; \; \forall \;\delta \hat{\bm{{\varepsilon}}}.
\end{equation}
Integrating by parts and using Eq. (\ref{eq_virtual-displacement}), Eq. (\ref{min_L2}) yields:
\begin{equation} \label{epsilon_cappello}
\hat{\bm{{\varepsilon}}}=\mathcal{G}^{-1} \bigg( \int_{\partial \mathit{\Omega}_E} 
(\mathbf{N}_E^\text{T}\mathbf{N}^P)^\text{T}
\mathbf{N}^V 
\text{ds}  \tilde{\mathbf{V}}-\int_{\mathit{\Omega}_E} 
(\mathbf{L}^\text{T}\mathbf{N}^P)^\text{T} \mathbf{v}_\mathit{h}
\text{dA} \bigg)
\end{equation}
where
\begin{equation} \label{eq_G_L2}
\mathcal{G} = 
\int_{\mathit{\Omega}_E}
(\mathbf{N}^P)^\text{T}
\mathbf{N}^P 
\text{dA},
\end{equation}
and $\mathbf{N}_E$ is the matrix containing the direction cosines of the outward unit normal vector on $\partial \mathit{\Omega}_E$.
The second integral in the r.h.s. of Eq. (\ref{epsilon_cappello}) collects the moments of the virtual displacement up to order $p-1$, and its evaluation would require $\mathbf{v}_\mathit{h}$ explicitly known in the interior of the element. To overcome this problem, in the VEM, the moments of the virtual displacement are assumed as internal degrees of freedom of the element, in addition to the external ones associated with the nodes on the element boundary.
Hence, denoting by $\hat{\mathbf{V}}$ the vector collecting the $p(p+1)$ local internal degrees of freedom (i.e. the moments of the virtual displacement up to order $p-1$: $\hat{V}_1=\int_{\mathit{\Omega}_E} v_{xh} \text{dA}$, $\hat{V}_2=\int_{\mathit{\Omega}_E} v_{yh} \text{dA}$, $\hat{V}_3=\int_{\mathit{\Omega}_E} x v_{xh} \text{dA}$, $\hat{V}_4=\int_{\mathit{\Omega}_E} y v_{xh} \text{dA}$,
$\hat{V}_5=\int_{\mathit{\Omega}_E} x v_{yh} \text{dA}$, ...),
 Eq. (\ref{epsilon_cappello}) can be rewritten as:
\begin{equation} \label{epsilon_cappello_2}
\hat{\bm{{\varepsilon}}}=\mathcal{G}^{-1}
( \tilde{\mathcal{B}} \tilde{\mathbf{V}} + \hat{\mathcal{B}}\hat{\mathbf{V}} ),
\end{equation}
where
\begin{equation} \label{eq_B_L2}
\tilde{\mathcal{B}} \tilde{\mathbf{V}} =
\int_{\partial \mathit{\Omega}_E} 
(\mathbf{N}_E^\text{T}\mathbf{N}^P)^\text{T}
\mathbf{N}^V 
\text{ds} \tilde{\mathbf{V}},
\qquad
\hat{\mathcal{B}} \hat{\mathbf{V}}=
-\int_{\mathit{\Omega}_E} 
(\mathbf{L}^\text{T}\mathbf{N}^P)^\text{T} \mathbf{v}_\mathit{h}
\text{dA}.
\end{equation}
As regards the second of Eq. (\ref{eq_B_L2}), it can be easily argued that, since $\hat{\mathbf{V}}$ collects the moments of the virtual displacement $\mathbf{v}_\mathit{h}$, $\hat{\mathcal{B}}$ is a sort of collocation operator collecting the scalar coefficients coming from the derivatives of the monomials in $\mathbf{N}^P$.

Summarizing, the total number of the element degrees of freedom is  $ n = 2mk + p (p + 1) $, 
 being $2mk$ the components of vector $\tilde{\mathbf{V}}$ and $p(p+1)$ those of $\hat{\mathbf{V}}$. Remembering that, in the standard VEM, $p=k-1$, it follows that internal degrees of freedom are present only if $k>1$.

Over the generic element, using Eqs. (\ref{eq_projector-operator}) and (\ref{epsilon_cappello_2}), the bilinear form in (\ref{PLV}) can be written in terms of projected strain as:
\begin{equation} \label{eq_discrete-bilinear}
\int_{\mathit{\Omega}_E}
[\bm{{\varepsilon}}^P(\delta \mathbf{v}_\mathit{h})]^\text{T}
\mathbf{C}
\bm{{\varepsilon}}^P(\mathbf{v}_\mathit{h})
\text{dA}
=
\int_{\mathit{\Omega}_E}
[\mathbf{N}^P\mathcal{G}^{-1}\mathcal{B} \; \delta \mathbf{V}]^\text{T}
\mathbf{C}
\mathbf{N}^P\mathcal{G}^{-1}\mathcal{B} \; \mathbf{V}
\text{dA},
\end{equation}
where
\begin{equation} \label{eq_B_TOTAL}
\mathcal{B} = 
\begin{bmatrix}
\tilde{\mathcal{B}} & \hat{\mathcal{B}}
\end{bmatrix},
\qquad
\mathbf{V} = 
\begin{bmatrix}
\tilde{\mathbf{V}} \\
 \hat{\mathbf{V}}
\end{bmatrix}.
\end{equation}
In  Eq. (\ref{eq_discrete-bilinear}), the consistent part of the stiffness matrix $\mathbf{K}_c$ can be recognized, which takes the form:
\begin{equation}\label{Kc}
\mathbf{K}_c=\mathcal{B}^{\text{T}}\mathcal{G}^{\text{-T}}\Big(\int_{\mathit{\Omega}_E} (\mathbf{N}^P)^\text{T}\mathbf{C}
\mathbf{N}^P 
\text{dA}\Big)\mathcal{G}^{-1}\mathcal{B}.
\end{equation}

\color{black}

\subsection{Stabilization term}\label{stabilization}
A suitable stabilizing term may be needed to preserve the coercivity of the system and avoid zero-energy modes. Indeed, the presence of a stability term is  standard for VEM \cite{beirao2013basic,beirao2013virtual,beirao2014hitchhiker,beirao2015virtual,artioli2017arbitrary}, even though stabilization may be not needed in some occasions, as specified in the following.
The formulation of the stabilization term herein presented is analogous to the one presented in \cite{artioli2017arbitrary}, to where the interested reader is referred.

Denoting by $s$ the stabilization order ($s=k$ in the standard displacement-based VEM \cite{artioli2017arbitrary}), 
the displacements field within the virtual element can be written both in terms of  vector-valued polynomials up to $s$ degree and also in terms 
of the virtual basis of the admissible displacements field of the element.
Accordingly, we call $\mathbf{D}$  the matrix associated to the change of basis
from the polynomial space
to the virtual functions space, i.e. the components of the matrix $\mathbf{D}$ are given by the evaluation of the polynomials on the degrees of freedom of the virtual element.

Consequently, the contribution of the stabilization in the stiffness matrix is
\begin{equation}\label{Ks}
\mathbf{K}_s=\tau \; \text{tr}(\mathbf{K}_c)\Big[
\mathbf{I}-\mathbf{D}\Big(\mathbf{D}^\text{T}\mathbf{D}\Big)^{-1}\mathbf{D}^\text{T}\Big],
\end{equation}
with $\tau$ positive real number (where the choice $\tau=1/2$ is common). It should be underlined that, as pointed out in \cite{artioli2017arbitrary}, the method is not much sensitive to the parameter $\tau$.
Finally, the trace term in (\ref{Ks}) guarantees the correct scaling of the energy with respect to the element size and material properties.

The total local stiffness matrix results in 
\begin{equation}\label{K}
\mathbf{K}=\mathbf{K}_c+\mathbf{K}_s.
\end{equation}

\subsection{Loading term}\label{loading}
Also for what concerns the loading term, direct reference to \cite{artioli2017arbitrary} is herein made, i.e. the loading term can be approximated by applying an integration rule based on vertexes for $k=1$, by using the average of the local load (through the first order moments) for $k=2$, and using higher order moments for $k>2$.

\subsection{Static condensation}\label{static}
Of course, in the same fashion of the finite element method, the standard Guyan condensation \cite{guyan1965reduction} for static problems (i.e. the static condensation) can be used to express the internal degrees of freedom in terms of those on the element boundary, so reducing the final degrees of freedom only to the latter. 

\section{Enhanced VEM formulation}\label{enahnced}
In this section, the possibility of modifying the standard VEM formulation presented in the previous section, operating on the hypotheses (i) and (ii) presented therein, is discussed.
In particular, with reference to the assumption (i) which provides for the use of the $L^2$ norm in Eq. (\ref {eq_requirement-projector}), it is observed that other choices are possible for this norm leading to different definitions of Eqs. (\ref{eq_G_L2})  and (\ref{eq_B_L2}).  
In particular, adopting an energy norm, the minimization of Eq. (\ref{eq_requirement-projector}) leads to: find $\bm{\varepsilon}^P \in \mathsf{P}_p(\Omega_E)$ such that
\begin{equation} \label{min_L2_b}
\int_{\mathit{\Omega}_E} [\bm{{\varepsilon}}^P-\bm{\varepsilon}(\mathbf{v}_\mathit{h})]^\text{T} \mathbf{C}\delta {\bm{\varepsilon}}^P \text{dA} = 0\; \; \; \forall \;\delta {\bm{\varepsilon}}^P\in \mathsf{P}_p(\Omega_E).
\end{equation}

\noindent Observing that $\mathbf{C}\bm{{\varepsilon}}^P$ is the stress in the element interior associated to the projected strain, it can be immediately observed that the  above equation expresses an energy equivalence condition.

Starting from Eq. (\ref{min_L2_b}) and following the same path outlined in the previous section, the following definitions are obtained:

\begin{equation} \label{eq_G_ENE_aa}
\mathcal{G} = 
\int_{\mathit{\Omega}_E}
({\mathbf{N}}^P)^\text{T}\mathbf{C}
{\mathbf{N}}^P 
\text{dA},
\end{equation}
\begin{equation} \label{eq_B_ENE_aa}
\tilde{\mathcal{B}} \tilde{\mathbf{V}} =
\int_{\partial \mathit{\Omega}_E} 
(\mathbf{N}_E^\text{T}\mathbf{C}{\mathbf{N}}^P)^\text{T}
\mathbf{N}^V 
\text{ds} 
\tilde{\mathbf{V}},
\qquad
\hat{\mathcal{B}} \hat{\mathbf{V}}=
-\int_{\mathit{\Omega}_E} 
(\mathbf{L}^\text{T}\mathbf{C}{\mathbf{N}}^P)^\text{T} \mathbf{v}_\mathit{h}
\text{dA}
.
\end{equation}

\noindent Even though it may appear a formal aspect, the use of the energy norm allows to highlight the mechanical meaning of the various quantities involved. In particular, noting that $\mathbf{C}{\mathbf{N}}^P$ can be interpreted as the representation basis of the stress in the element interior, it can be  argued that the second term of Eq. (\ref{eq_B_ENE_aa}) expresses for the work of the distributed volume forces in equilibrium with this stress field for the virtual displacement in the element.
Basing on this observation, and introducing  for later convenience:
\begin{equation} \label{NP_stress}
\breve{\mathbf{N}}^P=\mathbf{C}{\mathbf{N}}^P,
\end{equation}
\noindent the Eqs. (\ref{eq_G_ENE_aa}) and (\ref{eq_B_ENE_aa}) can be rewritten as:
\begin{equation} \label{eq_G_ENE}
\mathcal{G} = 
\int_{\mathit{\Omega}_E}
(\breve{\mathbf{N}}^P)^\text{T}\mathbf{C}^{-1}
\breve{\mathbf{N}}^P 
\text{dA},
\end{equation}
\begin{equation} \label{eq_B_ENE}
\tilde{\mathcal{B}} \tilde{\mathbf{V}} =
\int_{\partial \mathit{\Omega}_E} 
(\mathbf{N}_E^\text{T}\breve{\mathbf{N}}^P)^\text{T}
\mathbf{N}^V 
\text{ds} 
\tilde{\mathbf{V}},
\qquad
\hat{\mathcal{B}} \hat{\mathbf{V}}=
-\int_{\mathit{\Omega}_E} 
(\mathbf{L}^\text{T}\breve{\mathbf{N}}^P)^\text{T} \mathbf{v}_\mathit{h}
\text{dA}
.
\end{equation}

Accordingly, the consistent part of the stiffness matrix $\mathbf{K}_c$ in Eq. (\ref{Kc}) can be rewritten in terms of $\breve{\mathbf{N}}^P$, obtaining:
\begin{equation}\label{KcENERGY}
\mathbf{K}_c=\mathcal{B}^{\text{T}}\mathcal{G}^{\text{-T}}\Big(\int_{\mathit{\Omega}_E} (\breve{\mathbf{N}}^P)^\text{T}\mathbf{C}^{-1}
\breve{\mathbf{N}}^P 
\text{dA}\Big)\mathcal{G}^{-1}\mathcal{B}.
\end{equation}

As regards the hypothesis (ii) of the standard VEM, which provides that the polynomial degree $p$ in the representation basis of $\bm{{\varepsilon}}^P$ is linked to the order $k$ by the relationship (\ref{p_and_k}), i.e. $p = k-1$, some possibilities of enhancement of the standard VEM are presented in the following, focusing the attention directly on the definition of $\breve{\mathbf{N}}^P$ for convenience.


\color{black}

\subsection{Uncoupled polynomial representation}\label{div1}
One choice for the definition of $\breve{\mathbf{N}}^P$ consists in the adoption of an uncoupled polynomial representation of each component analogous to that of Eq. (\ref{eq_NP_div1}) but with the degree $p$ of the polynomial representation in $\breve{\mathbf{N}}^P$ not linked to the order $k$ by relationship (\ref{p_and_k}). 
This choice is indicated in the following as UnCoupled Polynomial representation or simply UCP.

Although various choices for the degree $p$ are possible, based on arguments similar to those for the optimal representations for stress and displacement in  mixed-stress finite elements
\cite{brezzi2012mixed,stolarski1987limitation},
 it can be  argued that an optimal choice for $p$ is such that the number of modes considered in $\breve{\mathbf{N}}^P$ fulfills the condition:
\begin{equation}\label{condition}
\mathrm{Number \: of  \: modes} \: 
\ge
n-3,
\end{equation}
i.e. when the number of modes considered in $\breve{\mathbf{N}}^P$ is greater than or equal to the element degrees of freedom  $n$ purged by the number of rigid motions (in this case 3).

As well as in the standard VEM, enriching $\breve{\mathbf{N}}^P$ up to a desired degree $p$ will lead to the origin of $p(p+1)$  moments of the virtual displacements field over the element (internal degrees of freedom), independently from the approximation degree of the element boundary $k$. Also in this case, a standard static condensation can be used to condensate these internal degrees of freedom.

Of course, the use of a projection on a higher order polynomial space  requires more integration points. Such aspect  increases the computational effort especially in nonlinear frameworks, where nonlinear constitutive laws must be evaluated at each integration point.

Finally, it is worth to note that the UCP enhancement could be also used in a standard  $L^2$ norm framework.
\color{black}


\subsection{Divergence-free polynomial representation}\label{div0}
Basing on the observations reported at the beginning of the Section, in case no distributed volume forces are present (i.e. $\mathbf{b}=\mathbf{0}$), the second term in Eq. (\ref{eq_B_ENE}) should be null. Accordingly, if $\mathbf{b}=\mathbf{0}$, a convenient choice for $\breve{\mathbf{N}}^P$ consists in the adoption of a matrix of preselected divergence-free modes that leads to:
\begin{equation} \label{eq_B_ENE_RED}
\hat{\mathcal{B}} \hat{\mathbf{V}}=
\mathbf{0}
.
\end{equation}
Such choice, 

which corresponds to a self-equilibrated representation of $\mathbf{C}\bm{{\varepsilon}}^P$ in the element interior,
\color{black}
is indicated in the following as Divergence-Free Polynomial representation or simply DFP. 
Particularly, the  components of the divergence-free $\breve{\mathbf{N}}^P$ matrix can be easily constructed by employing the well-established techniques usually adopted for hybrid stress finite elements.
In the plane elasticity framework, $\breve{\mathbf{N}}^P$ can be derived as 
\cite{benedetti2006posteriori}:
\begin{equation}  \label{eq_NP_div0}
\breve{\mathbf{N}}^P = 
\begin{bmatrix}
\breve{\mathbf{N}}^P_1\\
\breve{\mathbf{N}}^P_2\\
\breve{\mathbf{N}}^P_3
\end{bmatrix},
\;
\breve{\mathbf{N}}^P_1 = \frac{\partial^2\mathit{\phi}}{\partial y^2},
\;
\breve{\mathbf{N}}^P_2 = \frac{\partial^2\mathit{\phi}}{\partial x^2},
\;
\breve{\mathbf{N}}^P_3 = -\frac{\partial^2\mathit{\phi}}{\partial x \partial y}
\end{equation}
where the function $\mathit{\phi}$ is chosen as a complete polynomial function of degree $p+2$, being $p$ the desired degree of polynomial approximation considered in $\breve{\mathbf{N}}^P$.
Examples of $\breve{\mathbf{N}}^P$ with divergence-free modes are given in Appendix A.
As regards the optimal degree $p$ of the polynomial representation in $\breve{\mathbf{N}}^P$, an observation similar to that of the UPC case (Eq. (\ref{condition})) holds.

It is worth to note that, in the present DFP case, no internal degrees of freedom are needed to build the consistent part of the stiffness matrix of the element.
Therefore, the element degrees of freedom are   $n=2mk$.
 Accordingly, no static condensation is needed as no internal degrees of freedom origin in this case.

\paragraph{Remark 1}
For elements with order $k\ge2$ and in case of presence of distributed volume forces, a divergence-free $\breve{\mathbf{N}}^P$ would lead to loose the expected convergence rate of the standard $k$-th order VEM. Therefore, a special construction of the $\breve{\mathbf{N}}^P$ matrix could be performed in this case. Particularly, the representation is split into two parts: $(a)$ an uncoupled polynomial representation of each component complete up to degree $k-1$, $(b)$  divergence-free modes (following the aforementioned procedure) for higher-order terms up to the desired degree $p$. 
An example of $\breve{\mathbf{N}}^P$ for $k=2$ and $p=4$ is given in Appendix B.

Also in this case,  a standard static condensation procedure can be used to eliminate the internal degrees of freedom originated by the uncoupled polynomial part of degree $k-1$ in $\breve{\mathbf{N}}^P$. This choice is indicated in the following as HYbrid Polynomial representation or simply HYP.

\subsection{Some remarks about the stabilization term}
The enhancement of the Virtual Element Method presented in this section, beyond the enrichment of the polynomial representation of the strain in the element interior, has another interesting advantage. Indeed, when the condition (\ref{condition}) is fullfilled the method appears self-stabilized, i.e. no stabilization  is needed. In the standard VEM, this condition is only encountered in $k=1$ triangles. Conversely, this condition is fulfilled in several cases according to the enhanced VEM herein proposed, as shown in the following.
For example, an 8-node $k=1$ element ($m=8$) is self-stabilized if the degree of polynomial approximation is assumed  $p=3$, i.e. with 30 modes and  $n=28$ for UCP and with 18 modes and $n=16$ for DFP (for instance, it would not be self-stabilized with $p=2$, i.e. with 18 modes and  $n=22$ for UCP and with 12 modes and $n=16$ for DFP).

\section{Numerical results}
In this section, the  results of several numerical tests used to check the performance of the enhancements with respect to standard displacement-based VEM are shown and discussed.
In particular, the  strategy adopted for the VEM enhancement  is described in Section 5.1. The load cases (Load case A and Load case B) and meshes considered in the numerical tests are described in Section 5.2. Numerical results for Load case A and Load case B are shown and discussed in Section 5.3 and Section 5.4, respectively. Particularly, the enhanced VEM performance in nearly incompressible materials is discussed in Section 5.4.1.

\subsection{Adopted strategy}\label{strategy}
In this section, the strategy adopted for the VEM enhancement  is described.
Ultimately,  three enhancements (already introduced in Section 4) are considered in the following:
\begin{itemize}
	\item {UCP} $-$ uncoupled polynomial representation;
	\item {DFP} $-$ divergence-free polynomial representation;
	\item {HYP} $-$ hybrid polynomial representation.
\end{itemize}

For $k=1$, both UCP and DFP enhancements can be utilized. Particularly, static condensation of the $p(p+1)$  moments of the virtual displacements field is adopted in UCP.

For $k=2$, a distinction has to be made depending of the presence of a source term, i.e. null or non-null distributed forces. In the first case, UCP can be utilized in both cases and a static condensation of the $p(p+1)$ scalar moments is adopted. In the second case, DFP is used when distributed forces are null, otherwise HYP is used when distributed forces are non-null, following the discussion in Remark 1. Accordingly, no static condensation is implemented in DFP as no internal degrees of freedom origin in this case. Conversely, static condensation is used in HYP to condensate the internal degrees of freedom which origin by the uncoupled polynomial part of degree $k-1$.

The order $p$ of the polynomial representations is chosen as the minimum order which satisfies Eq. (\ref{condition}). 
Indeed, by further increasing the order $p$ beyond the minimum value which satisfies Eq. (\ref{condition}) no enhancements of the results are noted, as shown in Appendix C.
The values of  $p$ suggested for several enhanced VEM elements and used in the following are shown in Appendix D (Table \ref{table:strategy}), together with the suggested  stabilization terms.
In our discussion, we limit $p\le4$  as greater values of $p$ would introduce high-order polynomials which could increase the computational cost required to integrate Eqs. (\ref{eq_G_ENE}) and (\ref{eq_B_ENE})  and 
bring significant numerical errors.  
However, the limitation $p\le4$ does not appear significant for practical applications as this allows to have self-stabilized $k=2$ elements up to 48 degrees of freedom (Table \ref{table:strategy}).

In the following figures, ``VEM'' will indicate the standard displacement-based virtual element formulation. Also for the classical VEM, as well as for the enhanced one, static condensation is adopted.

\subsection{Load cases and considered meshes}
Two different load cases  are herein considered on a unit square domain $\Omega=(0,1)^2$, in plane stress conditions. The Young’s modulus and Poisson’s ratio of the material are set to be $E=2.5$ and $\nu=0.25$, respectively.
The following exact displacements field is assumed for the two load cases:
\begin{itemize}
	\item \underline{Load case A $-$ null distributed volume forces}
	\begin{equation}  \label{loadCaseA}
	\mathbf{v}(x,y)=
	\begin{bmatrix}
	-\frac{\displaystyle x^6}{80}+\frac{\displaystyle x^4y^2}{2}-\frac{13}{16}x^2y^4+\frac{3}{40}y^6 \\
	\frac{\displaystyle xy^5}{4}-\frac{5}{12}x^3y^3
	\end{bmatrix},
	\end{equation}
	leading to 
	$\mathbf{b}(x,y)=
	\begin{bmatrix}
	0 \\
	0
	\end{bmatrix}$.
	
	\item {\underline{Load case B $-$ non-null distributed volume forces}}
	\begin{equation}  \label{loadCaseB}
	\mathbf{v}(x,y)=
	\begin{bmatrix}
	x\sin(\pi x) \sin(\pi y) \\
	y\sin(\pi x) \sin(\pi y)
	\end{bmatrix},
	\end{equation}
	leading to 
	$\mathbf{b}(x,y)=
	\begin{bmatrix}
	\frac{11}{3}\pi^2 x\sin(\pi x) \sin(\pi y)-\frac{5}{3}\pi^2 y\cos(\pi x)\cos(\pi y)-7\pi\cos(\pi x)\sin(\pi y) \\
	\frac{11}{3}\pi^2 y\sin(\pi x) \sin(\pi y)-\frac{5}{3}\pi^2 x\cos(\pi x)\cos(\pi y)-7\pi\cos(\pi y)\sin(\pi x)
	\end{bmatrix}$.
\end{itemize}

%

\begin{figure}[!bth]
	\centering
	\subfigure[QUAD]{\label{fig:mesh1}\includegraphics[width=0.3\linewidth]{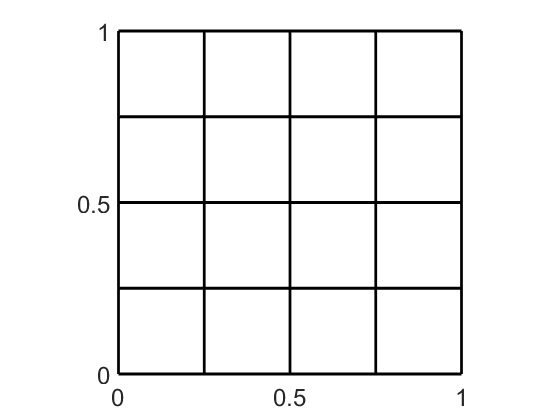}}
	\subfigure[RHOM]{\label{fig:mesh2}\includegraphics[width=0.3\linewidth]{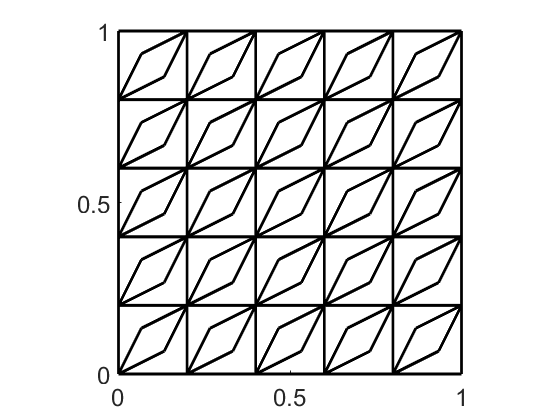}}\\
	\subfigure[HEXA]{\label{fig:mesh3}\includegraphics[width=0.3\linewidth]{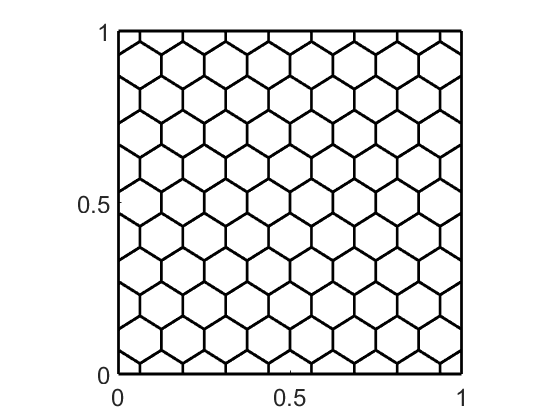}}
	\subfigure[WEBM]{\label{fig:mesh4}\includegraphics[width=0.3\linewidth]{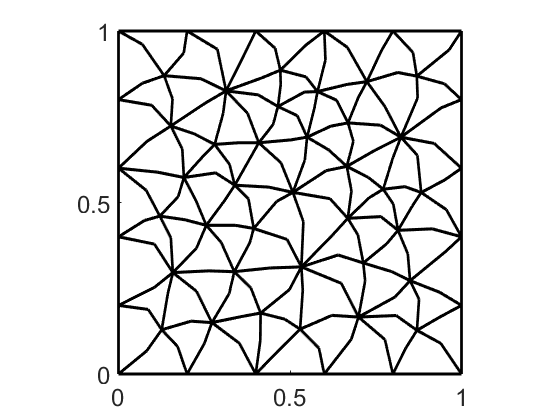}}
	\subfigure[DODE]{\label{fig:mesh5}\includegraphics[width=0.3\linewidth]{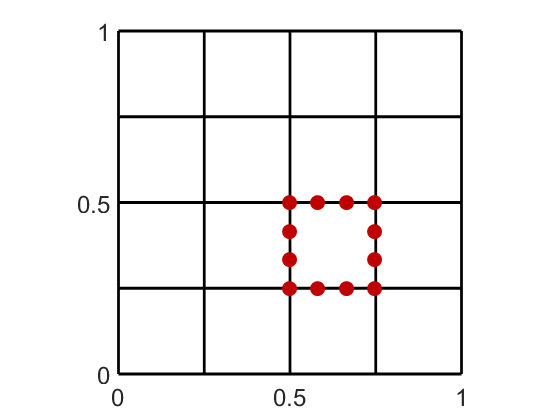}}
	\caption{Overview of considered meshes: (a) structured regular quadrilaterals, (b) convex/concave quadrilaterals, (c) structured regular hexagons, (d) convex/concave 6-vertex polygons, (e) structured 12-vertex polygons.}
	\label{fig:meshes}
\end{figure}

Five meshes have been considered to test the capability of the proposed enhanced VEM formulation (Figure \ref{fig:meshes}): structured regular quadrilaterals (QUAD, Figure \ref{fig:mesh1}), convex/concave quadrilaterals (RHOM, Figure \ref{fig:mesh2}), structured regular hexagons (HEXA, Figure \ref{fig:mesh3}),  convex/concave 6-vertex polygons (WEBM, Figure \ref{fig:mesh4}), structured 12-vertex polygons (DODE, Figure \ref{fig:mesh5}).

In the following graphs, ``log'' denotes logarithm with base 10, 
``dofs'' denotes,  for all the elements,  the degrees of freedom that result from the static condensation  (when adopted), while a classical energy norm is used to compute the error:

\begin{equation}\label{error}
||\mathrm{e}||=\frac{1}{||\bm{\varepsilon}^{EX}||}\left(\sum_{E} \left(\frac{1}{2} \int_{\Omega_E} (\bm{\varepsilon}^P-\bm{\varepsilon}^{EX})^T\mathbf{C}(\bm{\varepsilon}^P-\bm{\varepsilon}^{EX})\text{dA}\right)\right)^{1/2},
\end{equation}
where $\bm{\varepsilon}^{EX}$ is the exact strain field, 
$\bm{\varepsilon}^P$ is the projected strain, and 
$||\bm{\varepsilon}^{EX}||$ is the energy norm of the exact strain on the whole domain.

\color{black}
\clearpage

\subsection{Numerical results for Load case A}
\begin{figure}[!bth]
	\centering
	\subfigure[QUAD] {\label{fig:A_QUAD}\includegraphics[width=0.3\linewidth]{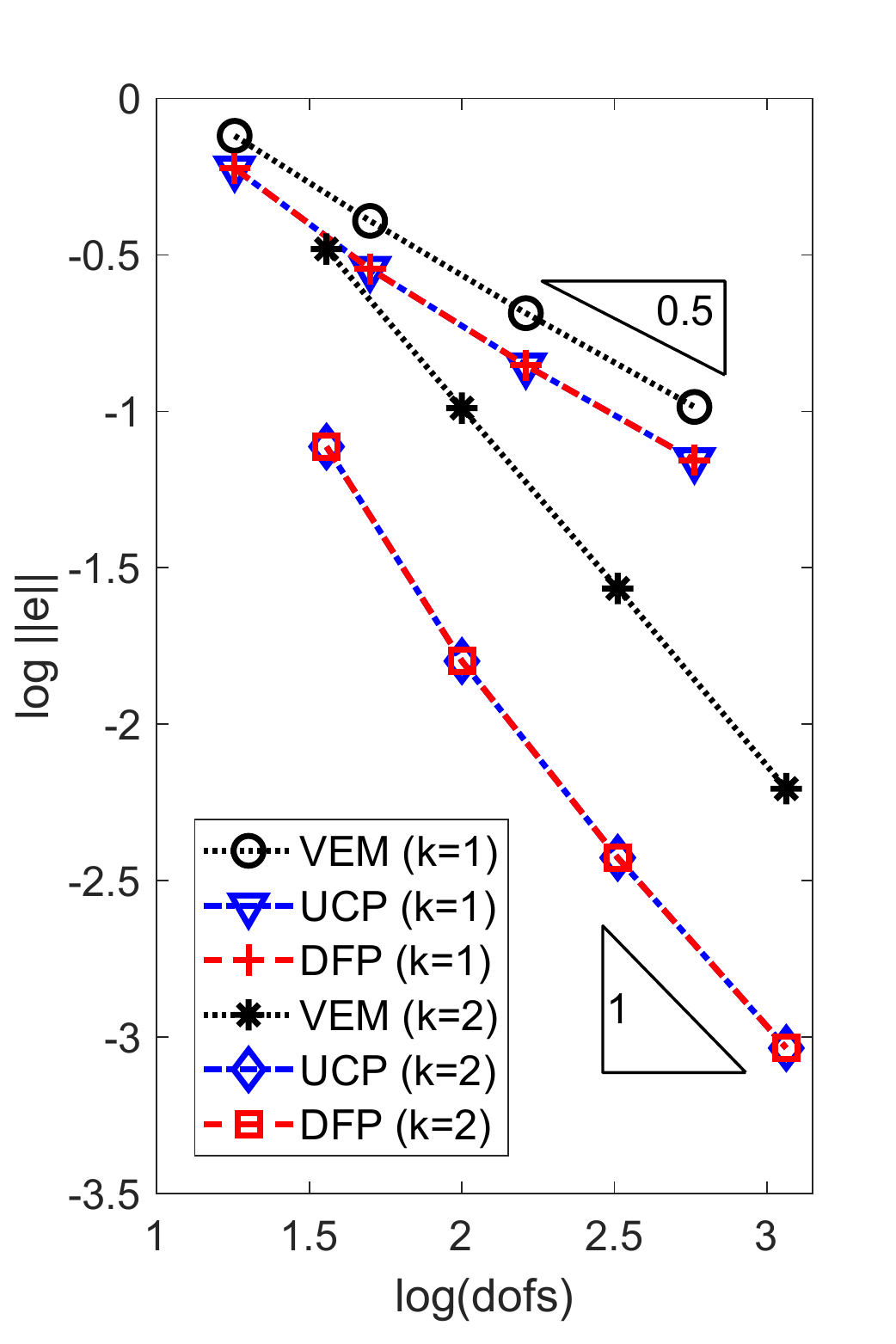}}
	\subfigure[RHOM] {\label{fig:A_RHOM}\includegraphics[width=0.3\linewidth]{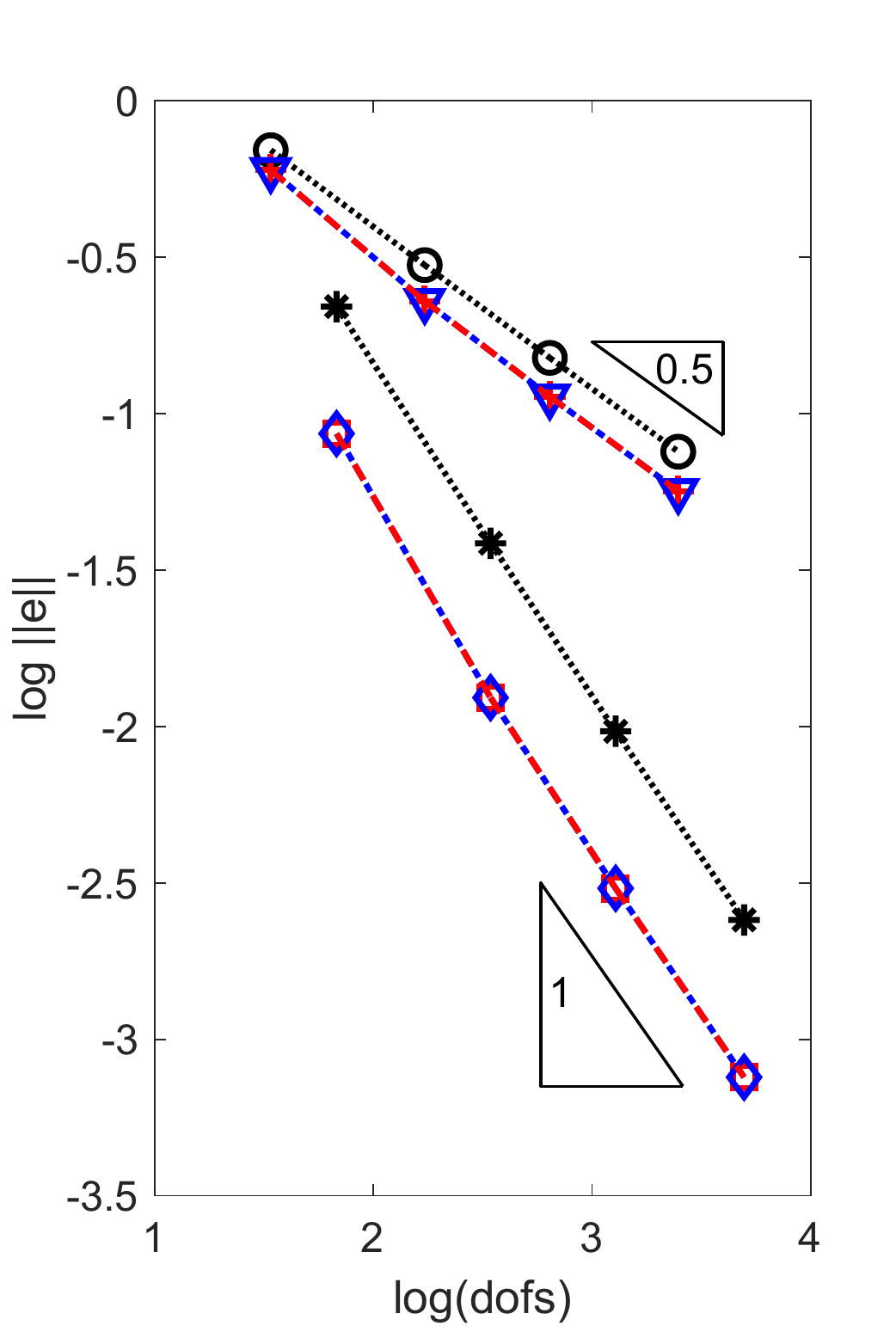}}\\
	\subfigure[HEXA] {\label{fig:A_HEXA}\includegraphics[width=0.3\linewidth]{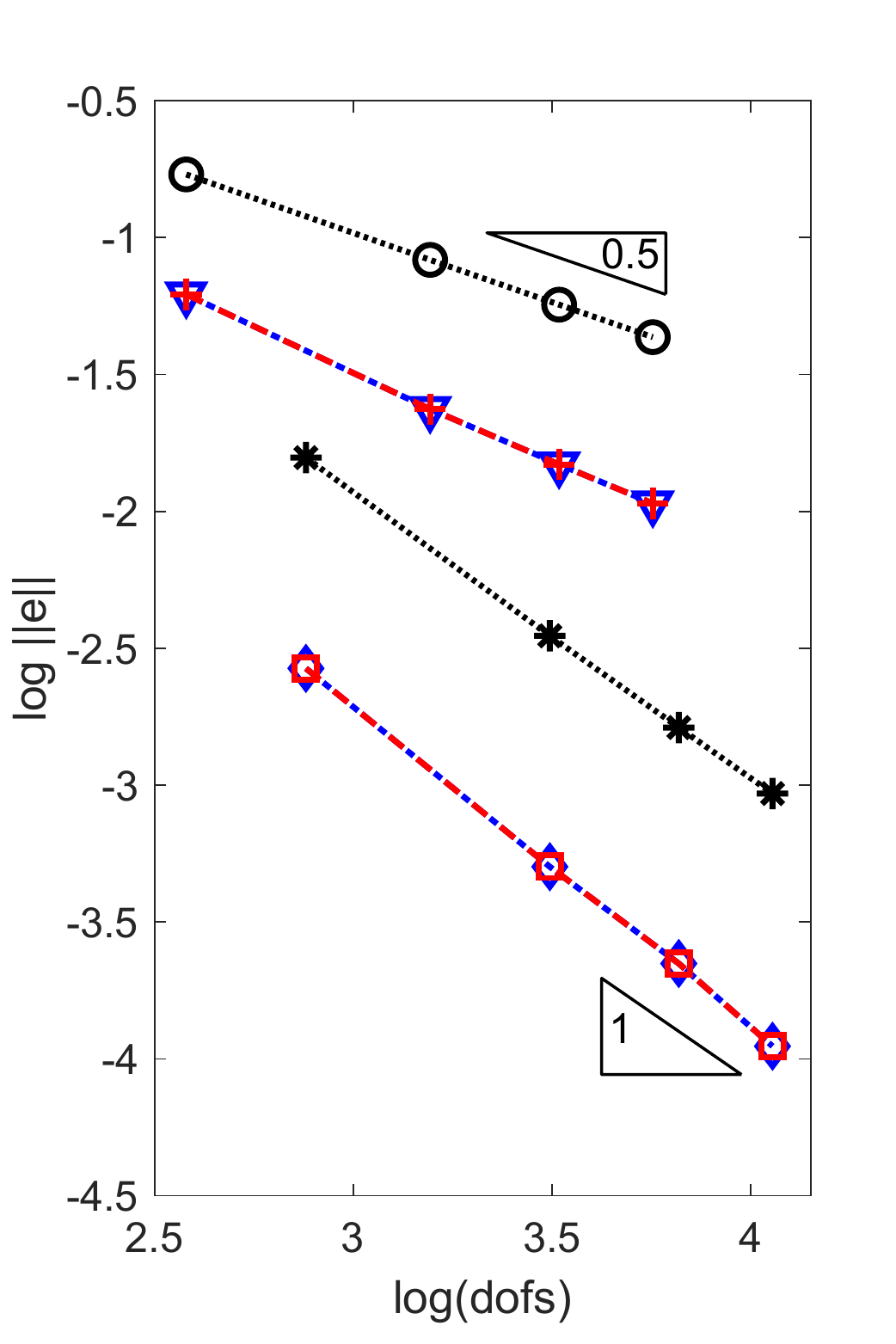}} 
	\subfigure[WEBM] {\label{fig:A_WEBM}\includegraphics[width=0.3\linewidth]{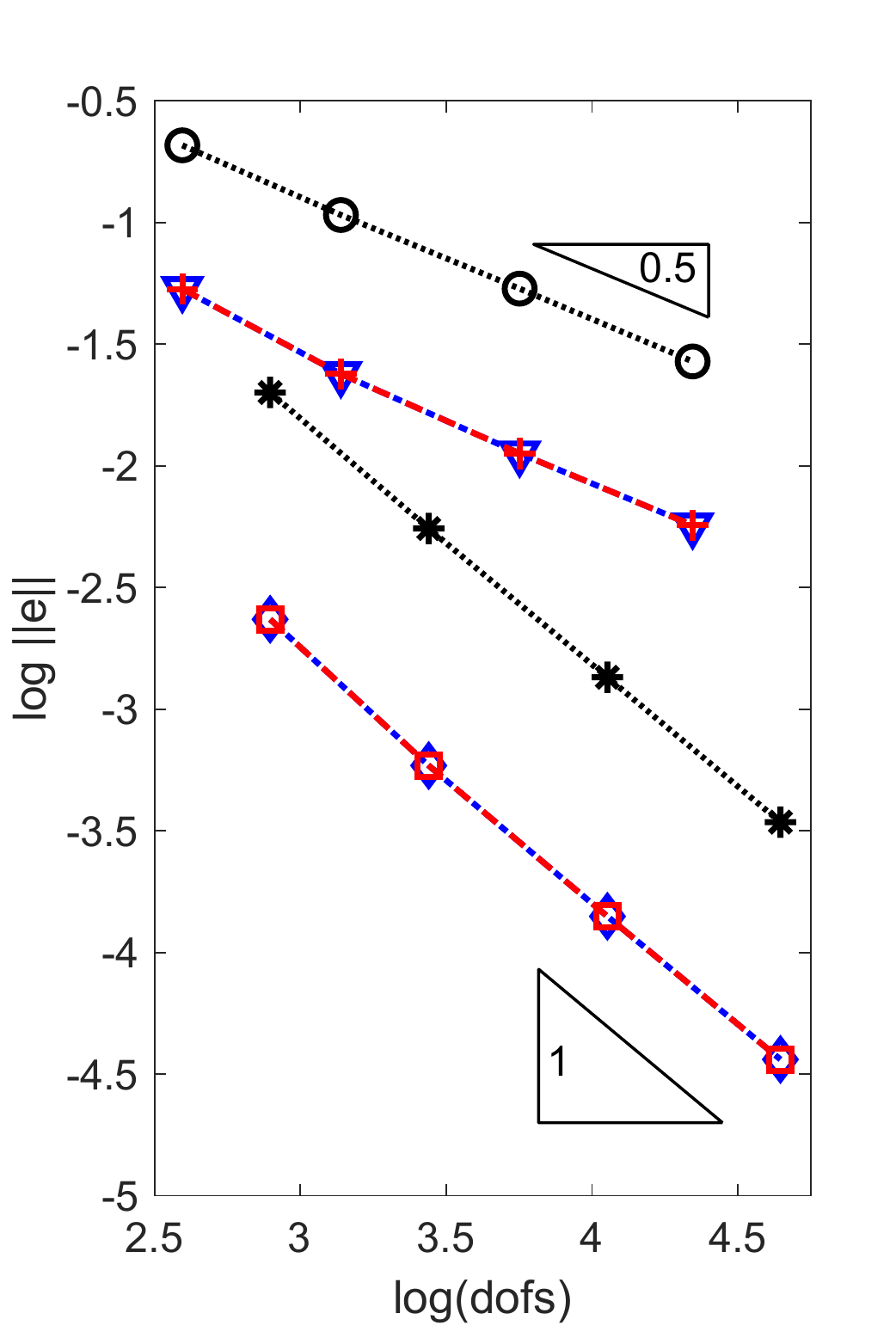}}
	\subfigure[DODE] {\label{fig:A_DODE}\includegraphics[width=0.3\linewidth]{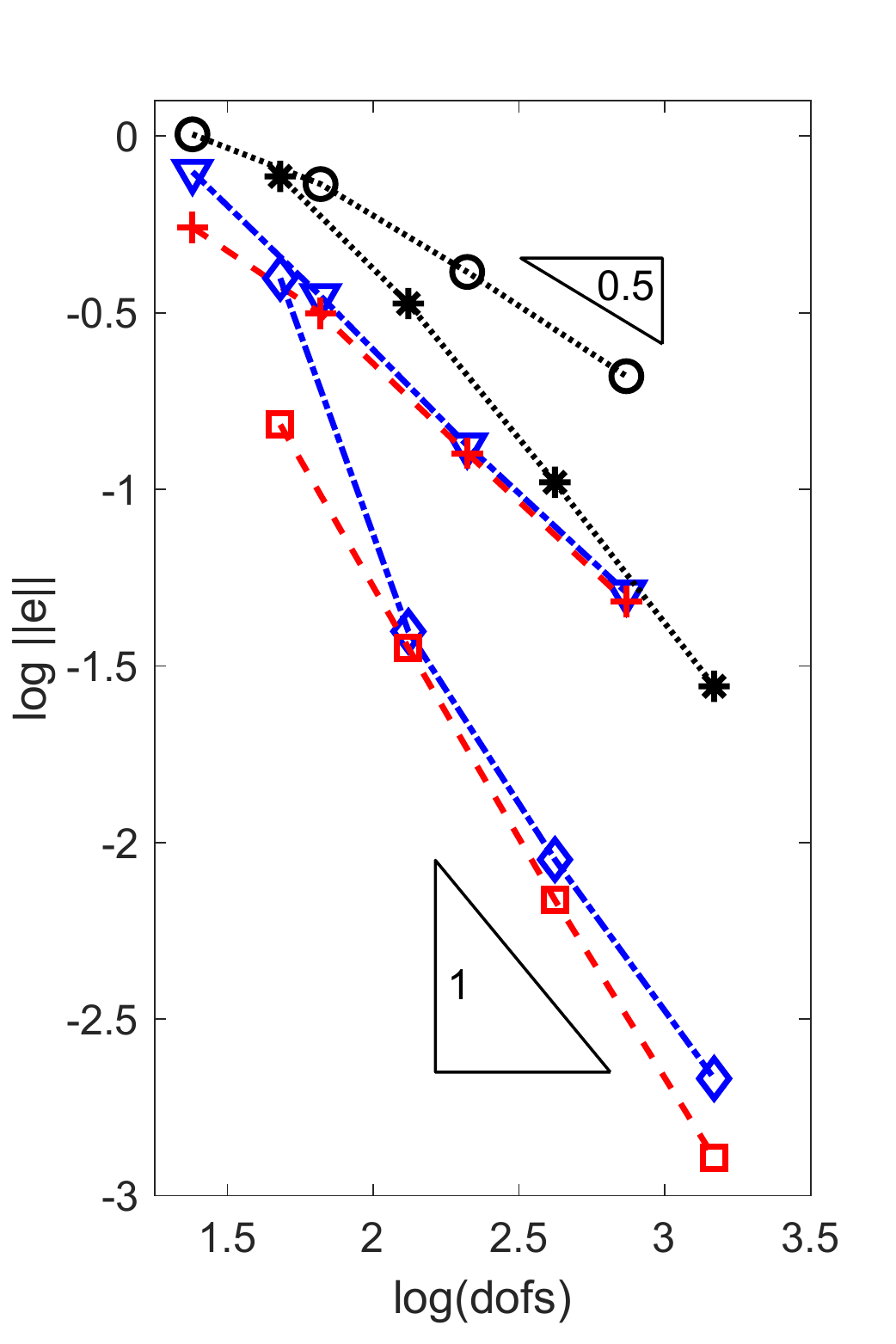}}
	\caption{Convergence for Load case A: (a) QUAD, (b) RHOM, (c) HEXA, (d) WEBM, and (e) DODE.}
	\label{fig:A_ALL}
\end{figure}

Concerning Load case A, which is characterized by null distributed forces, UCP and DFP enhancements have been considered for both $k=1$ and $k=2$ cases.
Convergence plots in terms of energy error (Eq. (\ref{error}))  are collected (for both $k=1$ and $k=2$ cases) in Figure \ref{fig:A_ALL} for QUAD (Figure \ref{fig:A_QUAD}), 
RHOM (Figure \ref{fig:A_RHOM}),
HEXA (Figure \ref{fig:A_HEXA}),
WEBM (Figure \ref{fig:A_WEBM}),
and  DODE (Figure \ref{fig:A_DODE}) meshes.

As it can be noted, for $k=1$   the results obtained with UCP and DFP enhancements show a  better performance, in terms of accuracy, than the standard displacement-based VEM for all meshes,
while keeping the optimal convergence rate of 0.5 for the energy error.
Also for $k=2$,  UCP and DFP enhancements show an increase of accuracy
 with respect to the standard  VEM, keeping the optimal convergence rate of 1 for  the energy error.
In Figure \ref{fig:A_ALL}, indeed,  the gain in accuracy increases by increasing the order $k$, i.e. passing from $k=1$ to $k=2$, for all meshes and for both UCP and DFP enhancements.
Also, the gain in accuracy appears to increase by increasing the number of vertexes $m$, for example passing from $m=4$ (e.g. QUAD, Figure \ref{fig:A_QUAD}) to $m=6$ (e.g. HEXA, Figure \ref{fig:A_HEXA}), for both $k=1$ and $k=2$ cases.

\color{black}
Particularly,  as could be expected due to the absence of distributed volume forces, UCP and DFP results coincide when stabilization is not utilized (i.e. in all cases except for DODE). Accordingly, DFP appears more efficient as it does not need static condensation.
When stabilization is utilized (e.g. in DODE, Figure \ref{fig:A_DODE}), the gain in accuracy obtained by UCP and DFP appears certainly comparable. In this case, results do not perfectly coincide as the stabilization in UCP and DFP  is different, as UCP and DFP are characterized by different dofs (note that static condensation is performed after stabilization, see Section \ref{stabilization}).

The gain in accuracy obtained with regular polygons appears  comparable with the one obtained with distorted polygons, e.g. compare QUAD (Figure \ref{fig:A_QUAD}) with RHOM (Figure \ref{fig:A_RHOM}) and HEXA (Figure \ref{fig:A_HEXA}) with WEBM (Figure \ref{fig:A_WEBM}). Therefore, the enhancement appears independent from mesh distortion.
Finally,
the DODE case (Figure \ref{fig:A_DODE}) highlights that both UCP and DFP show a significant enhancement in terms of accuracy with respect to standard VEM, even when stabilization is introduced (given the high number of vertexes $m=12$ and the will to keep $p\le4$ to prevent numerical round-off issues). Accordingly, the suggested stabilization terms in Table \ref{table:strategy} of Appendix D for polygons with high number of vertexes appear to not compromise the benefit of these enhancements.


\clearpage

\subsection{Numerical results for Load case B}

\begin{figure}[!bth]
	\centering
	\subfigure[QUAD] {\label{fig:B_QUAD}\includegraphics[width=0.3\linewidth]{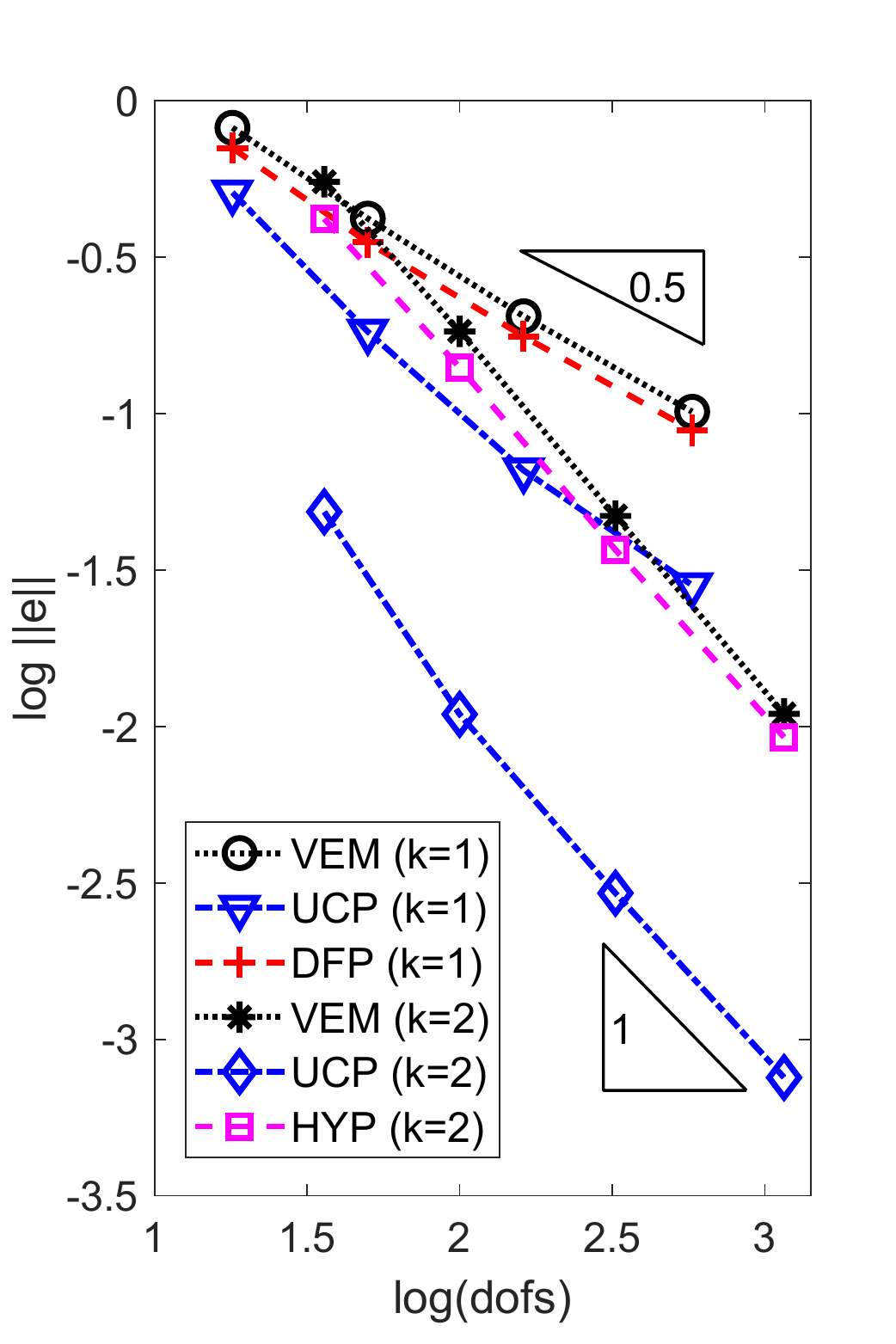}}
	\subfigure[RHOM] {\label{fig:B_RHOM}\includegraphics[width=0.3\linewidth]{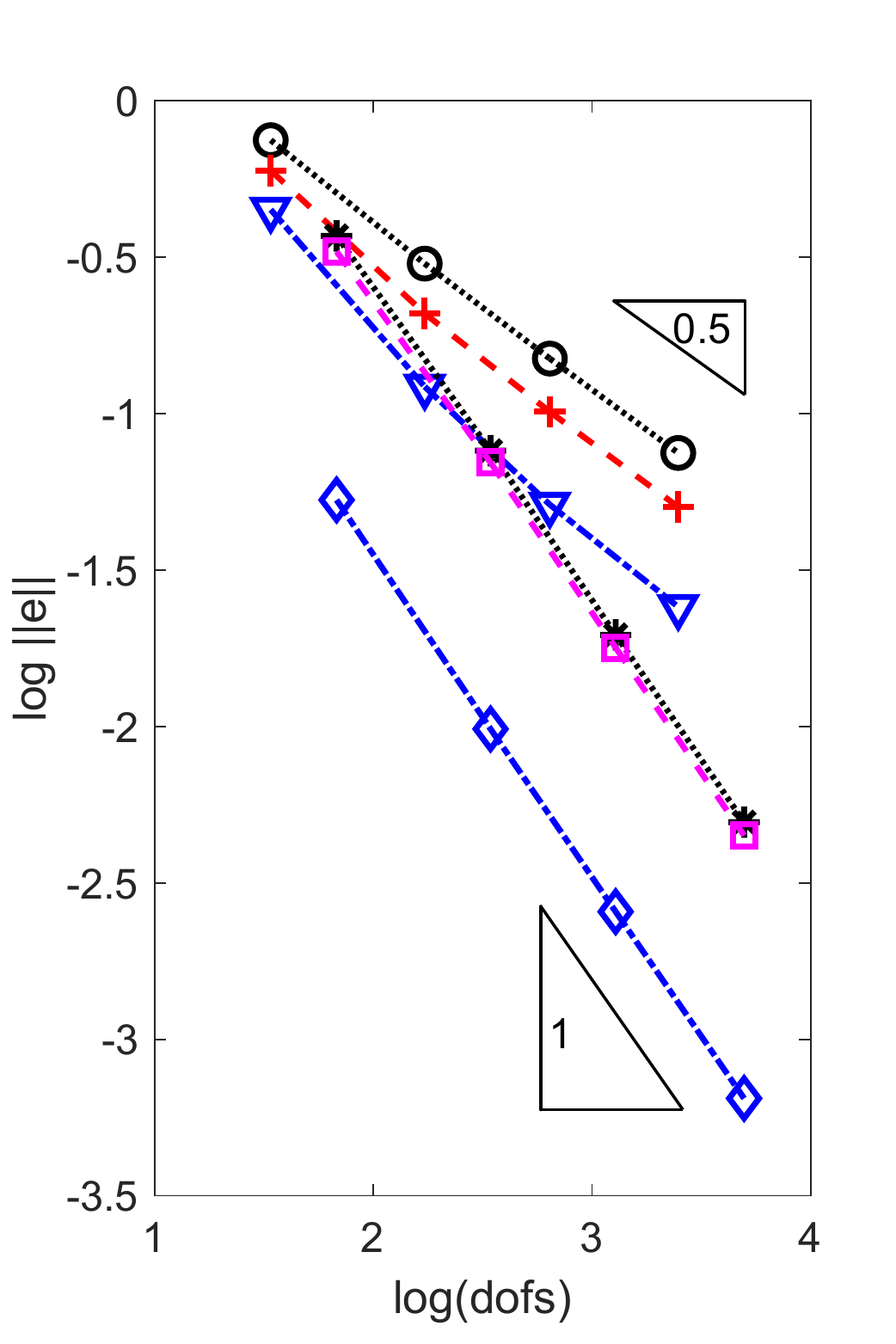}}\\
	\subfigure[HEXA] {\label{fig:B_HEXA}\includegraphics[width=0.3\linewidth]{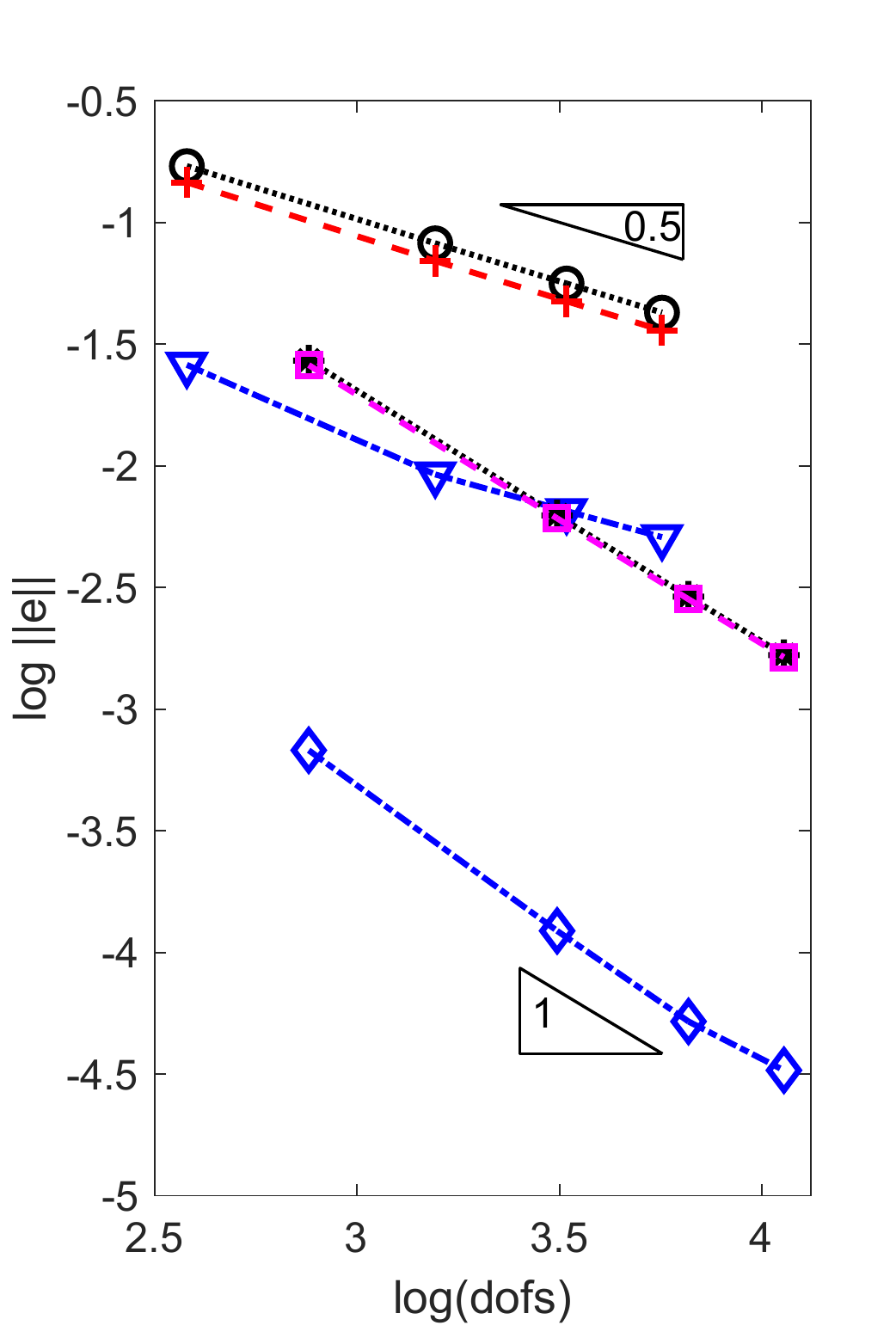}} 
	\subfigure[WEBM] {\label{fig:B_WEBM}\includegraphics[width=0.3\linewidth]{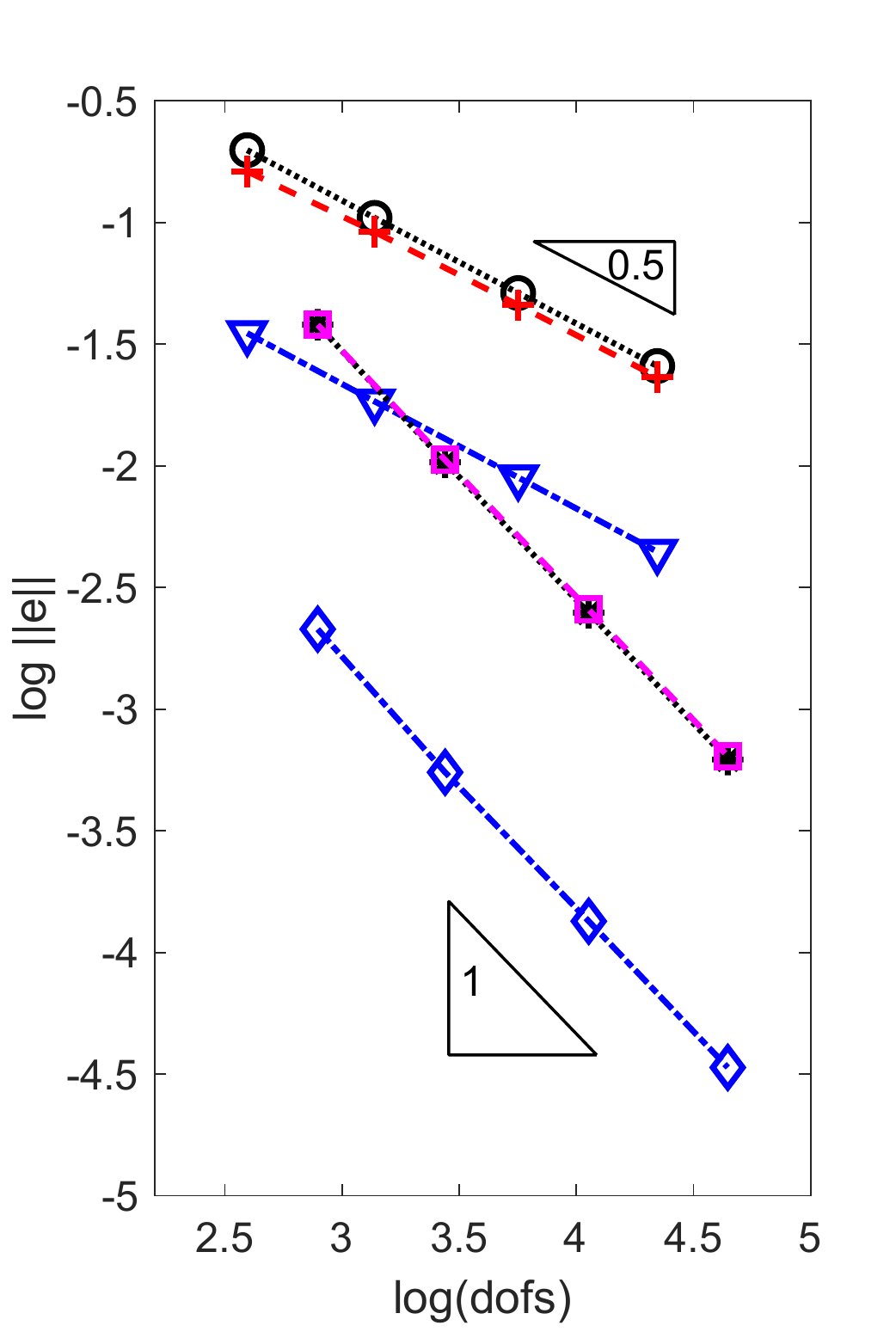}}
	\subfigure[DODE] {\label{fig:B_DODE}\includegraphics[width=0.3\linewidth]{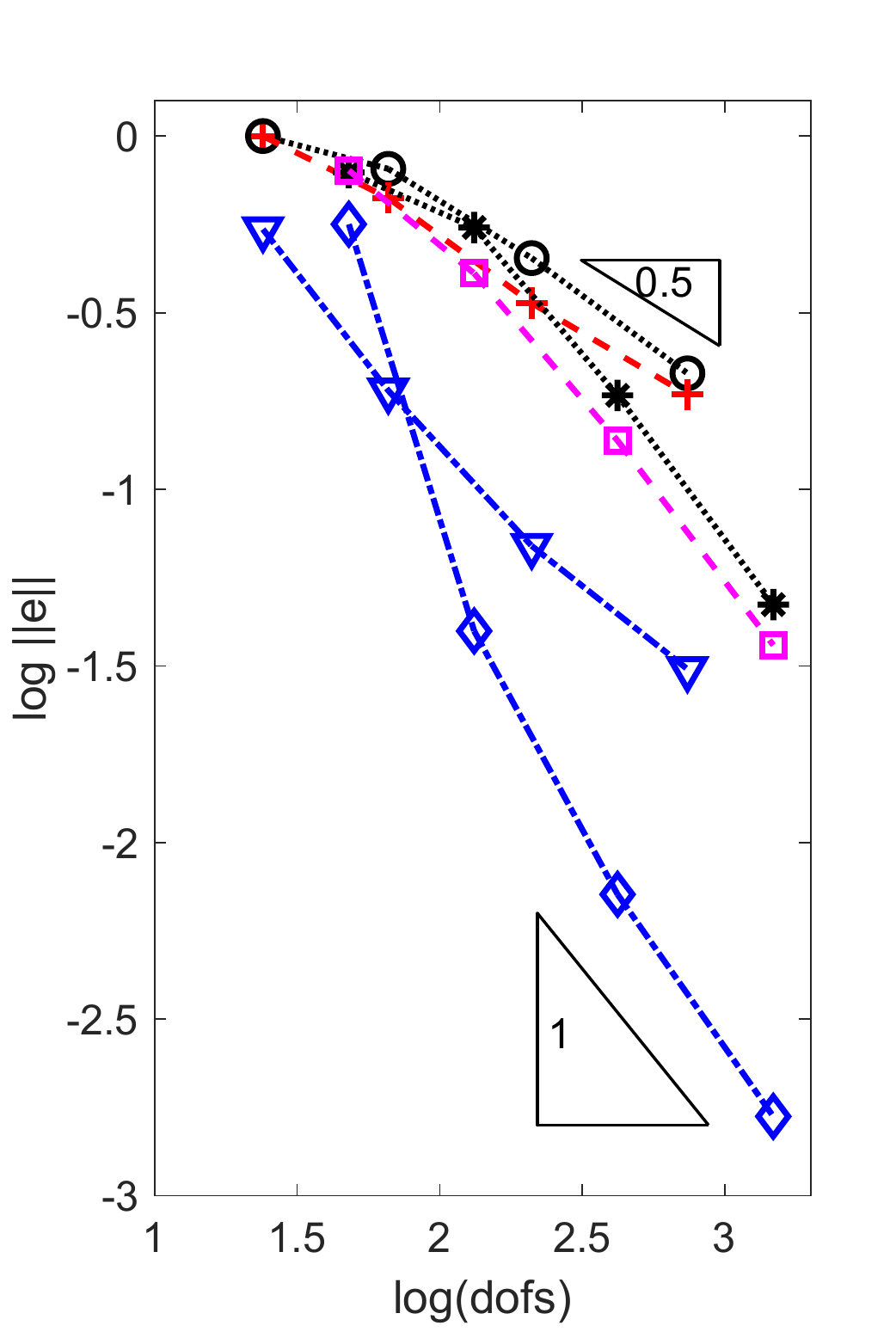}}
	\caption{Convergence for Load case B: (a) QUAD, (b) RHOM, (c) HEXA, (d) WEBM, and (e) DODE.}
	\label{fig:B_ALL}
\end{figure}

Concerning Load case B, which is characterized by non-null distributed forces, UCP and DFP enhancements have been considered for the $k=1$ case, while UCP and HYP enhancements have been considered for the $k=2$ case.
Convergence plots of energy error (Eq. (\ref{error}))  are collected (for both $k=1$ and $k=2$ cases) in Figure \ref{fig:B_ALL} for QUAD (Figure \ref{fig:B_QUAD}), 
RHOM (Figure \ref{fig:B_RHOM}),
HEXA (Figure \ref{fig:B_HEXA}),
WEBM (Figure \ref{fig:B_WEBM}),
and  DODE (Figure \ref{fig:B_DODE}) meshes, respectively.

The results obtained with UCP show a significant increase of accuracy with respect to the standard displacement-based VEM for all meshes, analogously to Load case A. Also in this case, the  optimal convergence rate is obtained with UCP for both $k=1$ and $k=2$ (Figure \ref{fig:B_ALL}), and the gain in accuracy increases by increasing the order $k$ for all meshes (Figure \ref{fig:B_ALL}).

\color{black}
As for Load case A, the gain in accuracy obtained in Load case B with UCP in regular polygons appears comparable with the one obtained in distorted polygons, e.g. compare QUAD (Figure \ref{fig:B_QUAD}) with RHOM (Figure \ref{fig:B_RHOM}) and HEXA (Figure \ref{fig:B_HEXA}) with WEBM (Figure \ref{fig:B_WEBM}). Also in this case, the enhancement appears independent from mesh distortion,  and the gain in accuracy increases by increasing the number of vertexes $m$ (in the cases without stabilization).

For the sake of comparison, 
 the convergence curves obtained with DFP (for the $k=1$ case) and HYP (for the $k=2$ case) enhancements
 are  reported
in Figure \ref{fig:B_ALL} as well. 
As can be noted, the DFP and HYP convergence curves keep the  optimal convergence rate of  the energy error for both $k=1$ and $k=2$. However, as could be expected due to the presence of non-null distributed volume load, their gain of accuracy with respect to the standard VEM solutions appears significantly more limited than the one observed for UCP. Even, the HYP solution coincides with the standard $k=2$ VEM in the HEXA (Figure \ref{fig:B_HEXA}) and WEBM  (Figure \ref{fig:B_WEBM}) cases.

Therefore,  the UCP enhancement represents an optimal solution in case of non-null distributed forces (i.e. in case of source terms), as it allows a considerable gain of accuracy, sensibly larger than the one observed in DFP and HYP. 

\clearpage
\subsubsection{Nearly incompressible material}
\begin{figure}[!bth]
	\centering
	\subfigure[QUAD $k=1$] {\label{fig:B_INCOMPR_QUAD}\includegraphics[width=0.350\linewidth]{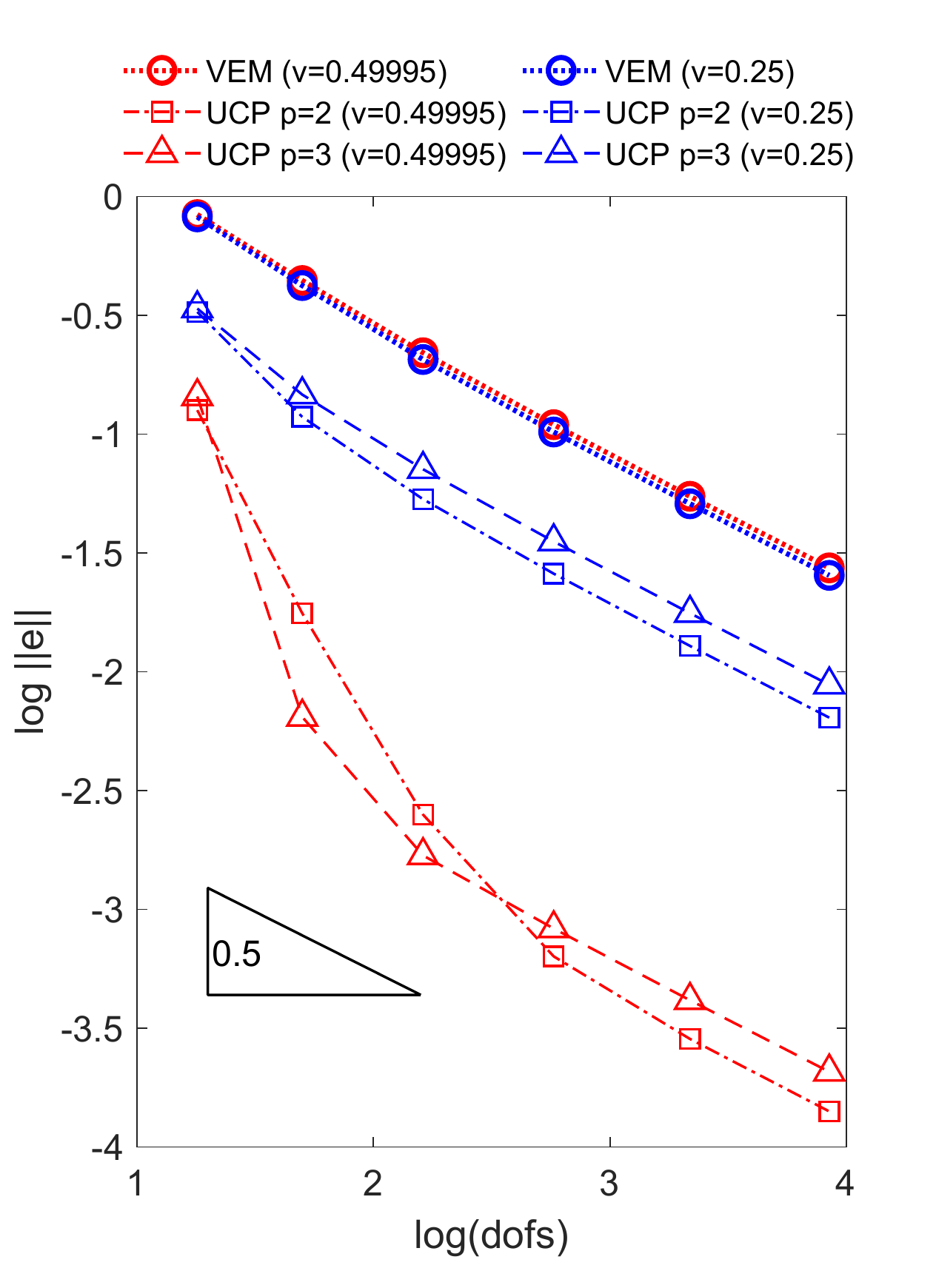}}
	\subfigure[HEXA $k=1$]
	{\label{fig:B_INCOMPR_HEXA}\includegraphics[width=0.350\linewidth]{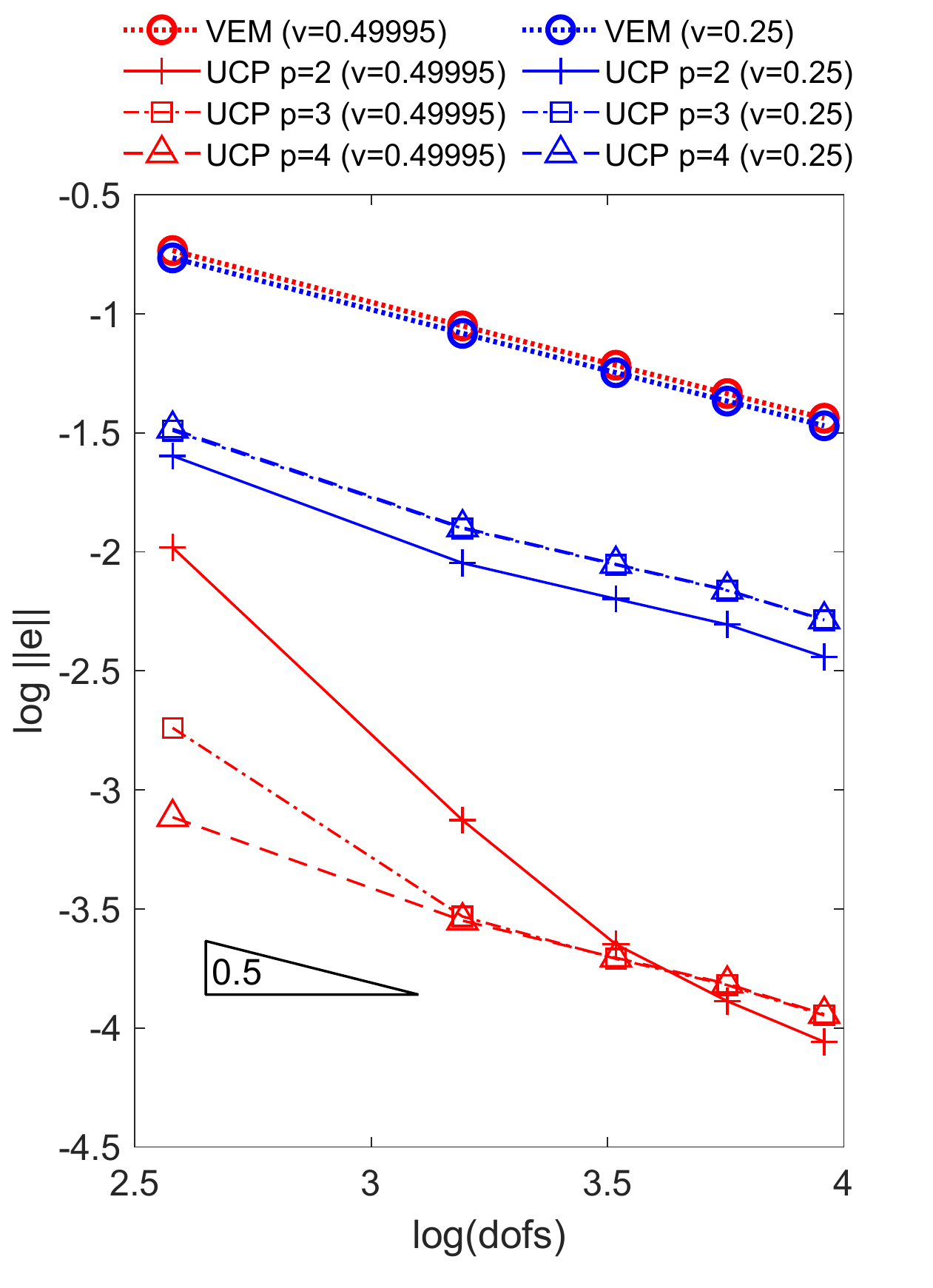}}
	\caption{Convergence for Load case B. Comparison between $\nu=0.25$ and $\nu=0.49995$ cases for different values of $p$: (a) QUAD $k=1$, and (b) HEXA $k=1$.}
	\label{fig:B_INCOMPR}
\end{figure}

Here,  UCP  is also tested in the case of a nearly incompressible material, by setting $\nu=0.49995$ (Figure \ref{fig:B_INCOMPR}) and plane strain conditions.
The comparison of  convergence curves  for Load case B for both $\nu=0.25$ and $\nu=0.49995$ cases, for different values of $p$,
are collected  in Figure \ref{fig:B_INCOMPR_QUAD} for QUAD $k=1$ and  in Figure \ref{fig:B_INCOMPR_HEXA} for HEXA $k=1$.

\begin{figure}[!bth]
	\centering
	\subfigure[VEM] {\label{fig:MAPS_VEM}\includegraphics[width=0.40\linewidth]{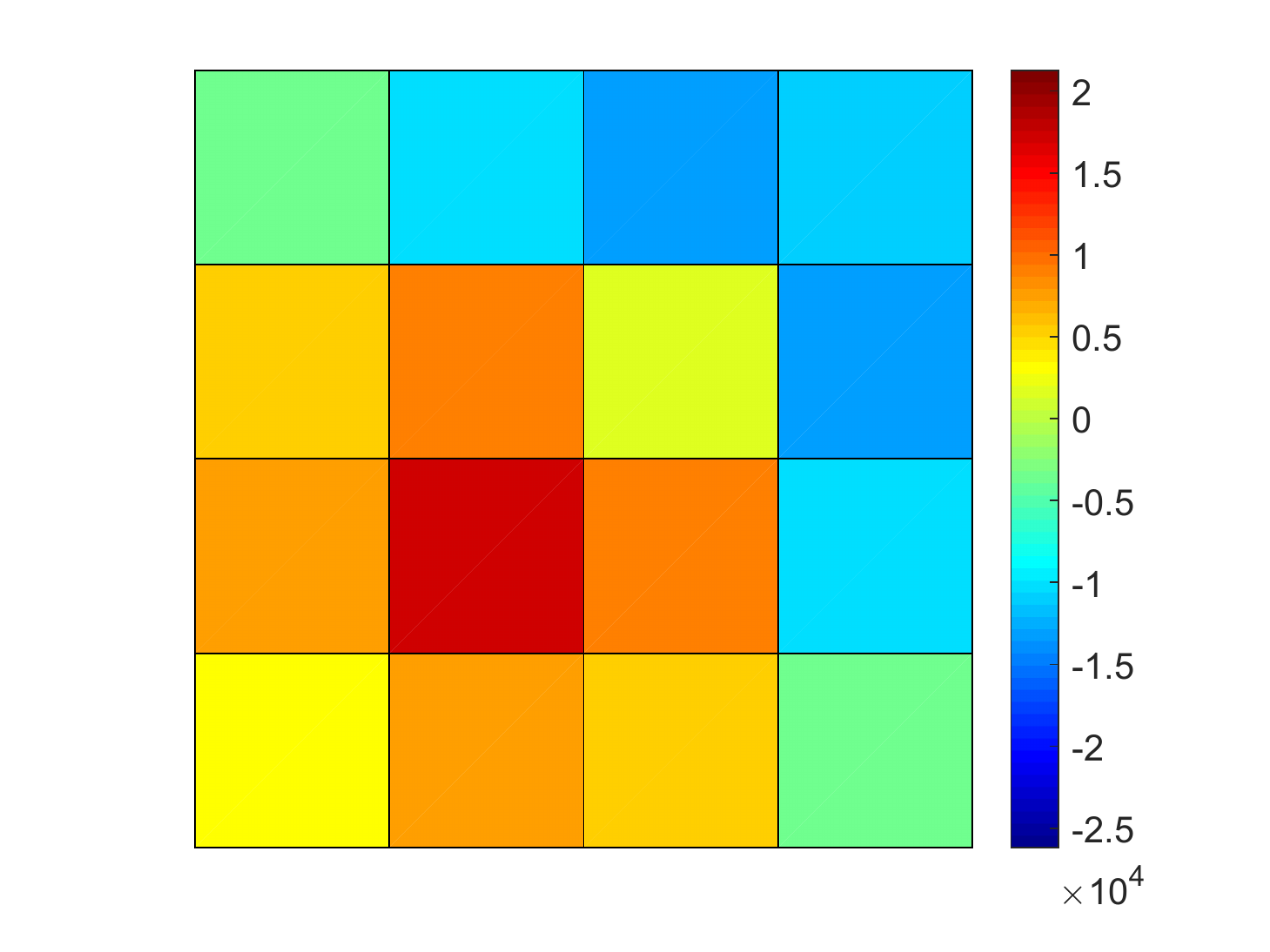}}
	\subfigure[UCP $p=2$] {\label{fig:MAPS_p2}\includegraphics[width=0.400\linewidth]{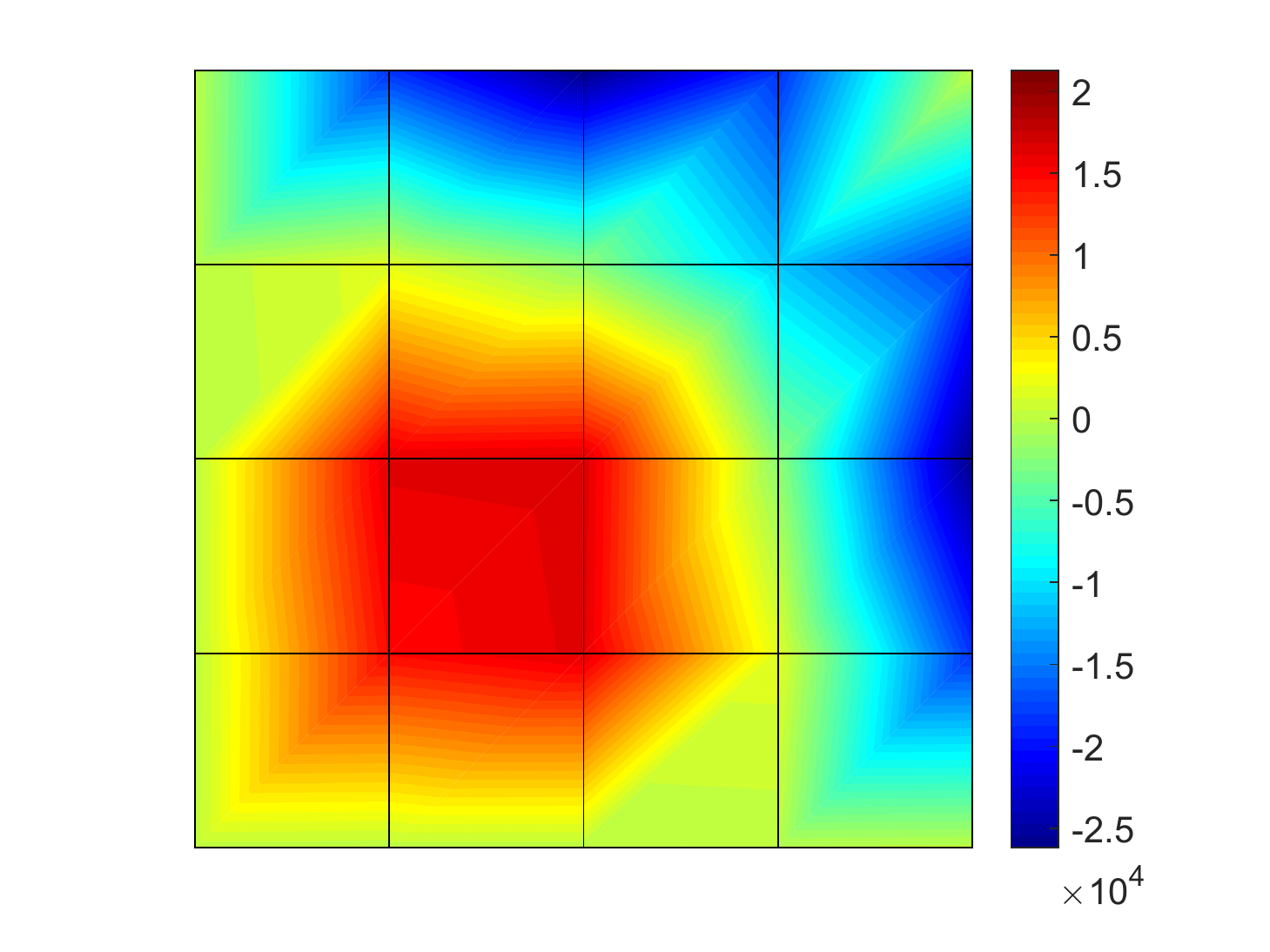}}
	\subfigure[UCP $p=3$] {\label{fig:MAPS_p3}\includegraphics[width=0.400\linewidth]{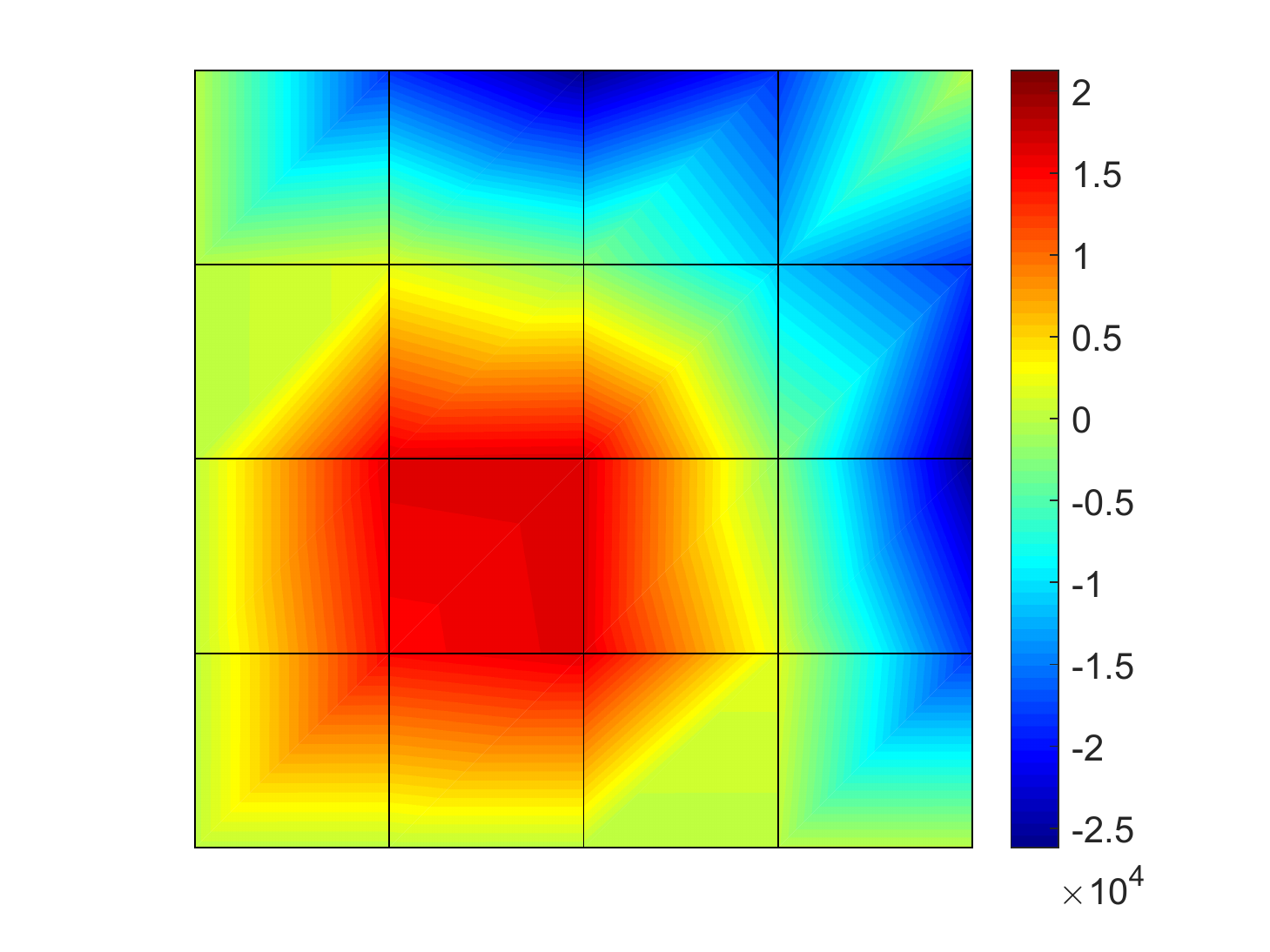}}
	\subfigure[Exact] {\label{fig:MAPS_exact}\includegraphics[width=0.400\linewidth]{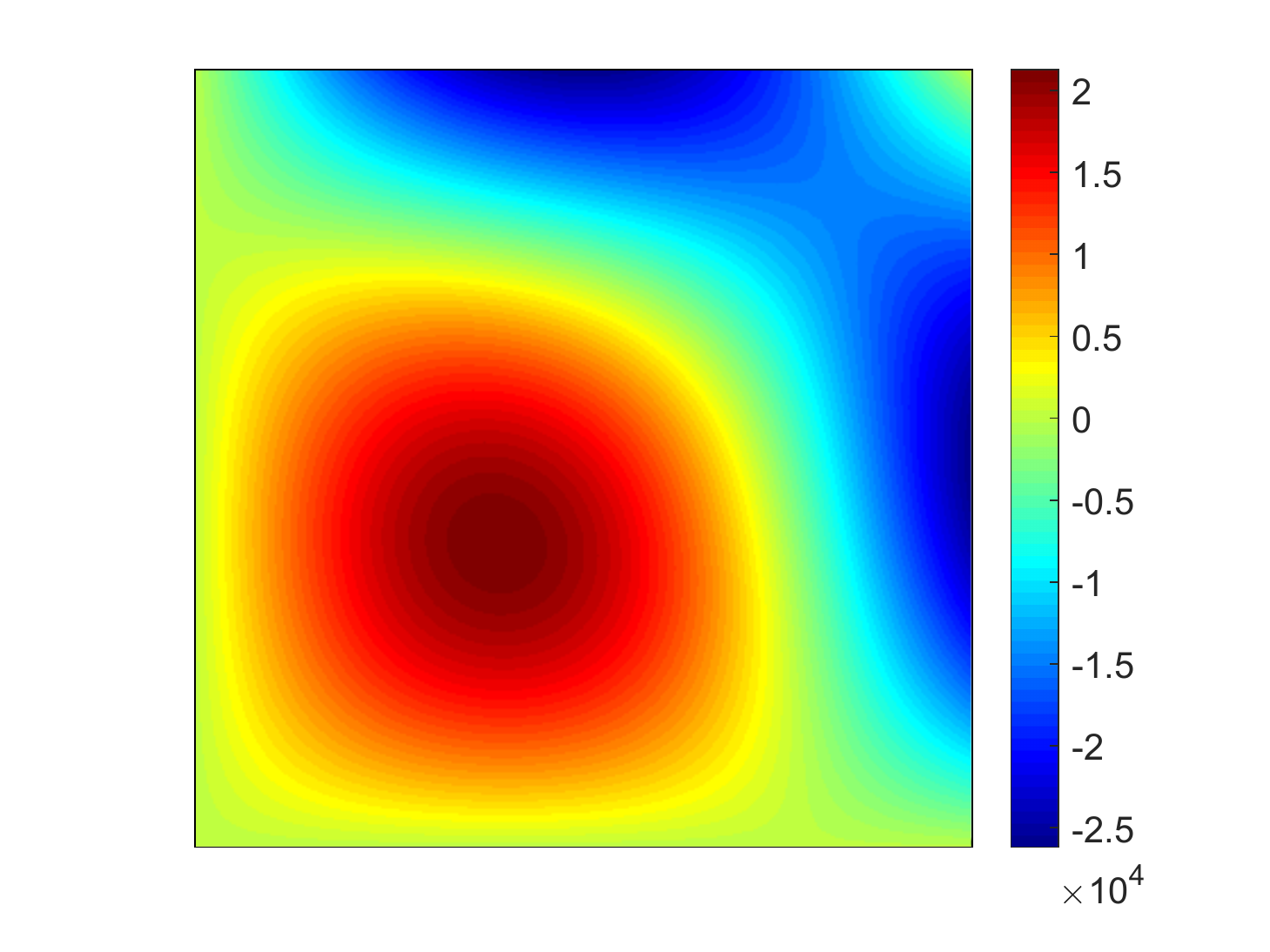}}
	\caption{Pressure maps for  Load case B with $\nu=0.49995$: QUAD $k=1$ (a) VEM,  (b) UCP $p=2$, (c) UCP $p=2$, and (d) exact solution.}
	\label{fig:MAPS}
\end{figure}

Although the enhancement of the VEM formulation presented in this paper has not been specifically developed for nearly incompressible materials, it  appears  that this enhancement can be effectively used also in this case (Figure \ref{fig:B_INCOMPR}).
Indeed, thanks to the higher order internal representation, the increase of accuracy obtained in the case of nearly incompressible material ($\nu=0.49995$)  appears remarkable for both QUAD (Figure \ref{fig:B_INCOMPR_QUAD}) and HEXA (Figure \ref{fig:B_INCOMPR_HEXA}) meshes (in this case even greater than the $\nu=0.25$ one). This increase of accuracy appears even clearer in Figure \ref{fig:MAPS}, where the pressure maps  for Load Case B with $\nu=0.49995$ are shown for a $4\times4$ QUAD $k=1$ mesh: standard VEM (Figure \ref{fig:MAPS_VEM}), UCP $p=2$ (Figure \ref{fig:MAPS_p2}), UCP $p=3$ (Figure \ref{fig:MAPS_p3}), as well as the exact solution  (Figure \ref{fig:MAPS_exact}). As can be noted, enhanced solutions show considerably better accuracy than standard VEM which predicts, as expected, constant pressure within each element.

\color{black}

\clearpage
\section{Conclusions}\label{conclusions}

In this paper, an enhanced VEM formulation has been proposed for plane elasticity.
In the proposed formulation, elements have enhanced strain representation within the element, while keeping the same degree of the displacement approximation functions on the element boundary of the standard VEM formulation.

The enhanced VEM formulation has been proposed with both uncoupled (UCP) and divergence-free (DFP) polynomial representations, which have been tested through several numerical examples. On the one hand, natural serendipity elements generated by the DFP enhancement represented an optimal solution in case of null distributed forces. On the other hand, the UCP enhancement  always showed a significant gain of accuracy with respect to standard VEM also in presence of distributed forces.

It appears worth to highlight here that the UCP enhancement can be also used in the standard $L^2$ norm framework. Conversely, the DFP enhancement requires the use of an energy norm framework to obtain optimal results  in case of null distributed forces.

\color{black}
To conclude, numerical results  showed the capability of the enhanced VEM formulation 
 to (i) considerably increase accuracy (with respect to standard VEM) while keeping the optimal convergence rate, (ii) bypass the need of stabilization terms in many practical cases, (iii) obtain natural serendipity elements  in many practical cases, and (vi) effectively treat also the case of nearly incompressible materials.

%

%
%
%
%
\section*{Appendix A}\label{appA}
In the following, examples of $\breve{\mathbf{N}}^P$ with divergence-free polynomial representation (DFP) are given for:
\begin{itemize}
	\item $p=1$
	\begin{equation}  \label{eq_NP_div0p1}
	\breve{\mathbf{N}}^P = 
	\begin{bmatrix}
	1& 0& 0& y& 0& x&   0\\
	0& 1& 0& 0& x& 0&   y\\
	0& 0& 1& 0& 0& -y& -x
	\end{bmatrix};
	\end{equation}
	
	\item $p=2$
	\begin{equation}  \label{eq_NP_div0p2}
	\breve{\mathbf{N}}^P = 
	\begin{bmatrix}
	1& 0& 0& y& 0& x&   0& 2xy&    0&     x^2& y^2&   0\\
	0& 1& 0& 0& x& 0&   y &   0& 2xy&   y^2&   0& x^2\\
	0& 0& 1& 0& 0& -y& -x&  -y^2& -x^2& -2xy&   0&   0
	\end{bmatrix};
	\end{equation}
	
	\item $p=3$
	
	\begin{equation}  \label{eq_NP_div0p3}
	\tiny
	\breve{\mathbf{N}}^P = 
	\begin{bmatrix}
	1& 0& 0& y& 0& x&   0&     x^2& y^2&   0&       0& -2xy&       x^3& y^3&   0& 3x^2y&         0& -3y^2x\\
	0& 1& 0& 0& x& 0&   y&     y^2&   0& x^2& -2xy&       0&  3y^2x&   0& x^3&       y^3& -3x^2y&         0\\
	0& 0& 1& 0& 0& -y& -x& -2xy&   0&   0&     x^2&     y^2& -3x^2y&   0&   0& -3y^2x&       x^3&       y^3
	\end{bmatrix}; 
	\small
	\end{equation}

	\item $p=4$
	
	\begin{equation}  \label{eq_NP_div0p4}
	\tiny
	\breve{\mathbf{N}}^P = 
	\left[
	\begin{smallmatrix}
	1& 0& 0& y& 0& x&   0&     x^2& y^2&   0&       0& -2xy&       x^3& y^3&   0& 3x^2y&         0& -3y^2x& y^4& 0&-4y^3x& 0& 2x^3y& x^4&6x^2y^2\\
	0& 1& 0& 0& x& 0&   y&     y^2&   0& x^2& -2xy&       0&  3y^2x&   0& x^3&       y^3& -3x^2y&         0&0&x^4& 0& -4x^3y&2y^3x&6x^2y^2&  y^4\\
	0& 0& 1& 0& 0& -y& -x& -2xy&   0&   0&     x^2&     y^2& -3x^2y&   0&   0& -3y^2x&       x^3&       y^3&0&0& y^4&x^4&-3x^2y^2&-4x^3y&-4y^3x
	\end{smallmatrix}	\right]  .
	\small
	\end{equation} 	
\end{itemize}

\section*{Appendix B}\label{appB}
Here, examples of $\breve{\mathbf{N}}^P$ are given for the case HYP with:
\begin{itemize}
	
	\item $p=3$
	\begin{equation}  \label{eq_NP_div01p3}
	\breve{\mathbf{N}}^P = 
	\left[
	\begin{smallmatrix}
	1&0&0&x&0&0&y&0&0&     x^2& y^2&   0&       0& -2xy&       x^3& y^3&   0& 3x^2y&         0& -3y^2x\\
	0& 1&0&0&x&0&0&y&0&     y^2&   0& x^2& -2xy&       0&  3y^2x&   0& x^3&       y^3& -3x^2y&         0\\
	0& 0& 1&0&0&x&0&0&y& -2xy&   0&   0&     x^2&     y^2& -3x^2y&   0&   0& -3y^2x&       x^3&       y^3
	\end{smallmatrix}	\right]  ;
	\small
	\end{equation} 	
	
	\item $p=4$
	\begin{equation}  \label{eq_NP_div01p4}
	\tiny
	\breve{\mathbf{N}}^P = 
	\left[
	\begin{smallmatrix}
	1&0&0&x&0&0&y&0&0&     x^2& y^2&   0&       0& -2xy&       x^3& y^3&   0& 3x^2y&         0& -3y^2x& y^4& 0&-4y^3x& 0& 2x^3y& x^4&6x^2y^2\\
	0& 1&0&0&x&0&0&y&0&     y^2&   0& x^2& -2xy&       0&  3y^2x&   0& x^3&       y^3& -3x^2y&         0&0&x^4& 0& -4x^3y&2y^3x&6x^2y^2&  y^4\\
	0& 0& 1&0&0&x&0&0&y& -2xy&   0&   0&     x^2&     y^2& -3x^2y&   0&   0& -3y^2x&       x^3&       y^3&0&0& y^4&x^4&-3x^2y^2&-4x^3y&-4y^3x
	\end{smallmatrix}	\right]  .
	\small
	\end{equation} 	
\end{itemize}

\section*{Appendix C}\label{appC}
\begin{figure}[!bth]
	\centering
	\subfigure[HEXA] {\label{fig:INF_HEXA1}\includegraphics[width=0.30\linewidth]{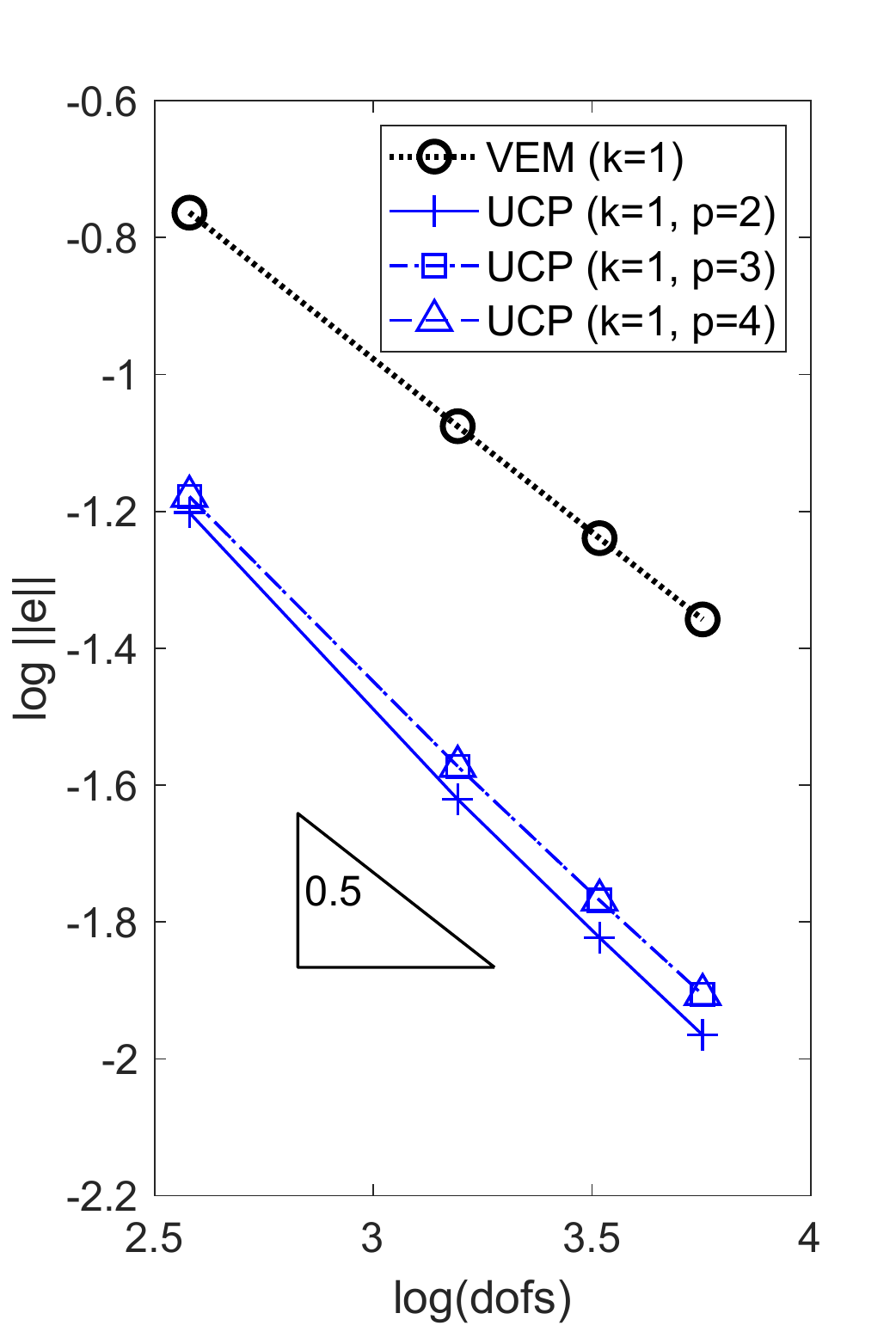}}
	\subfigure[QUAD]
	{\label{fig:INF_QUAD2}\includegraphics[width=0.30\linewidth]{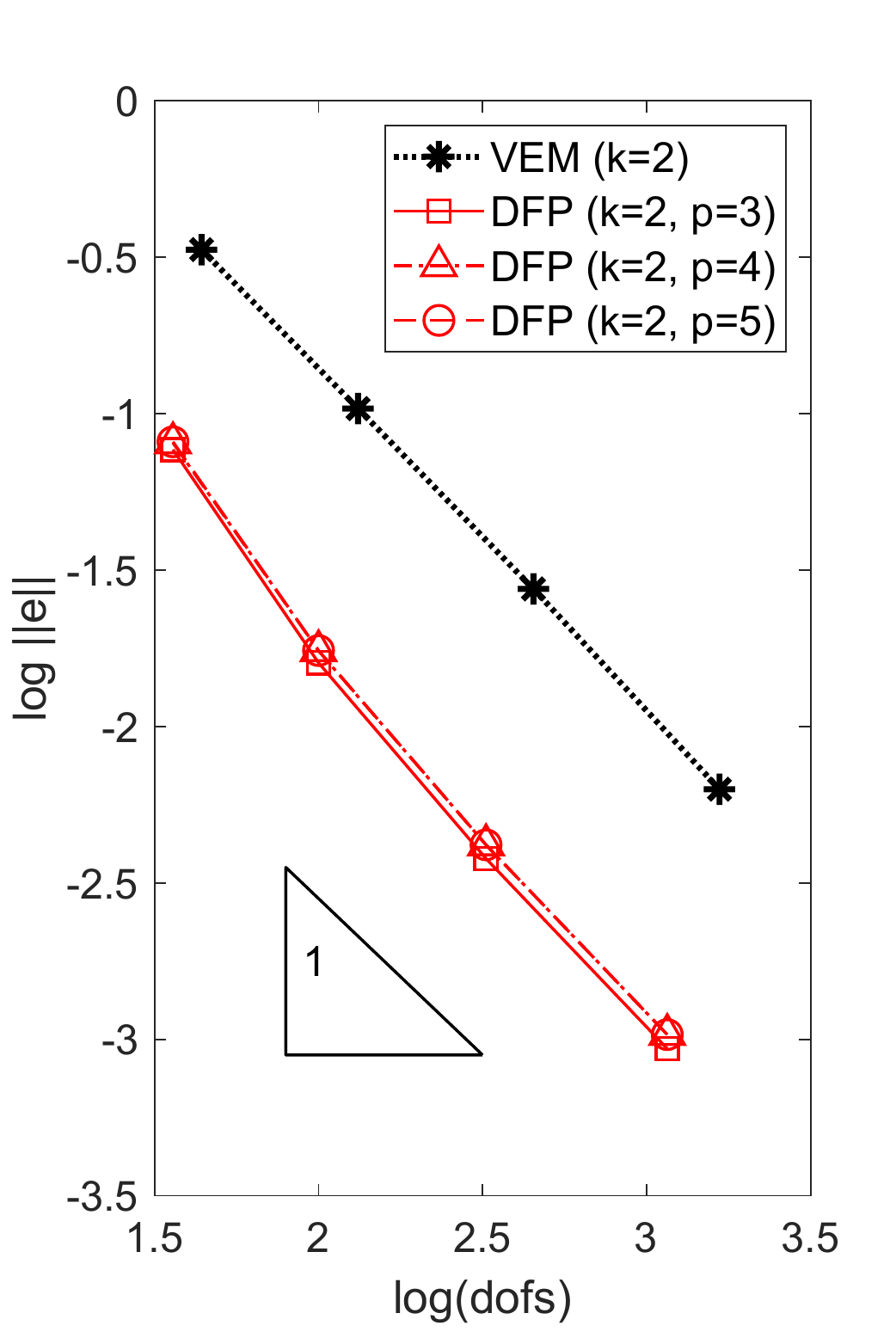}}
	\caption{Influence of $p$ on: (a) UPC (HEXA $k=1$), and (b) DFP (QUAD $k=2$), Load case A.}
	\label{fig:parametric}
\end{figure}

Figure \ref{fig:parametric} shows the influence of the adopted $p$ (in any case taken to satisfy Eq.  (\ref{condition})) on UPC (Figure\ref{fig:INF_HEXA1}) and DFP enhancements (Figure\ref{fig:INF_QUAD2}).

\section*{Appendix D}\label{appD}
Table \ref{table:strategy} shows the suggested values of $p$  for  various elements types.
\begin{table}[h!]
	\centering
	\caption{Values of $p$ suggested for  various elements. When Eq. (\ref{condition}) is satisfied, no stabilization is needed.
		$n$ denotes the total number of degrees of freedom including, when present, internal degrees of freedom. $s$ denotes the stabilization order.
		For the sake of comparison, in standard VEM once $k$ is adopted, it is assumed $p=k-1$ and $s=k$. }
	\label{table:strategy}
	\footnotesize
	\begin{tabular}{ |c|c|c|c|c|c|c| } %
		\hline
		$m$ (n. vertexes) & $k$ & $p$&Enhancement&$n$&$\breve{\mathbf{N}}^P$ modes& Stabilization \\
		\hline
		\multirow{5}{*}{4} & \multirow{2}{*}{1}  & \multirow{2}{*}{1} & UCP &10 & 9 & NO \\ 
		\cline{4-7}
		&   & & DFP &8 & 7 & NO \\ 
		\cline{2-7}
		& \multirow{3}{*}{2}  & \multirow{3}{*}{3} & UCP &28 & 30 & NO \\
		\cline{4-7}
		&   & & DFP &16 & 18 & NO \\ 
		\cline{4-7}
		&   & & HYP &18 & 20 & NO \\ 
		\hline
		\multirow{5}{*}{5} & \multirow{2}{*}{1}  & \multirow{2}{*}{1} & UCP &12 & 9 & NO \\ 
		\cline{4-7}
		&   & & DFP &10 & 7 & NO \\ 
		\cline{2-7}
		& \multirow{3}{*}{2}  & \multirow{3}{*}{3} & UCP &32& 30 & NO \\
		\cline{4-7}
		&   & & DFP &20 & 18 & NO \\ 
		\cline{4-7}
		&   & & HYP &22 & 20 & NO \\ 
		\hline
		\multirow{5}{*}{6} & \multirow{2}{*}{1}  & \multirow{2}{*}{2} & UCP &18 & 18 & NO \\ 
		\cline{4-7}
		&   & & DFP &12 & 12 & NO \\ 
		\cline{2-7}
		& \multirow{3}{*}{2}  & \multirow{3}{*}{4} & UCP &44& 45 & NO \\
		\cline{4-7}
		&   & & DFP &24 & 25 & NO \\ 
		\cline{4-7}
		&   & & HYP &26 & 27 & NO \\ 
		\hline
		\multirow{5}{*}{7} & \multirow{2}{*}{1}  & \multirow{2}{*}{2} & UCP &20 & 18 & NO \\ 
		\cline{4-7}
		&   & & DFP &14 & 12 & NO \\ 
		\cline{2-7}
		& \multirow{3}{*}{2}  & \multirow{3}{*}{4} & UCP &48& 45 & NO \\
		\cline{4-7}
		&   & & DFP &28 & 25 & NO \\ 
		\cline{4-7}
		&   & & HYP &30 & 27 & NO \\ 
		\hline
		\multirow{5}{*}{8} & \multirow{2}{*}{1}  & \multirow{2}{*}{3} & UCP &28 & 30 & NO \\ 
		\cline{4-7}
		&   & & DFP &16 & 18 & NO \\ 
		\cline{2-7}
		& \multirow{3}{*}{2}  & \multirow{3}{*}{4} & UCP &52& 45 & $s=4$ \\
		\cline{4-7}
		&   & & DFP &32 & 25 & $s=3$\\ 
		\cline{4-7}
		&   & & HYP &34 & 27 & $s=3$ \\ 
		\hline
		\multirow{5}{*}{9} & \multirow{2}{*}{1}  & \multirow{2}{*}{3} & UCP &30 & 30 & NO \\ 
		\cline{4-7}
		&   & & DFP &18 & 18 & NO \\ 
		\cline{2-7}
		& \multirow{3}{*}{2}  & \multirow{3}{*}{4} & UCP &56& 45 & $s=4$ \\
		\cline{4-7}
		&   & & DFP &36 & 25 & $s=3$\\ 
		\cline{4-7}
		&   & & HYP &38 & 27 & $s=3$ \\ 
		\hline
		\multirow{5}{*}{10} & \multirow{2}{*}{1}  & \multirow{2}{*}{3} & UCP &32 & 30 & NO \\ 
		\cline{4-7}
		&   & & DFP &20 & 18 & NO \\ 
		\cline{2-7}
		& \multirow{3}{*}{2}  & \multirow{3}{*}{4} & UCP &60& 45 & $s=4$ \\
		\cline{4-7}
		&   & & DFP &40 & 25 & $s=3$\\ 
		\cline{4-7}
		&   & & HYP &42 & 27 & $s=3$ \\ 
		\hline
		\multirow{5}{*}{$m>10$} & \multirow{2}{*}{1}  & \multirow{2}{*}{3} & UCP &$2m+12$ & 30 & $s=4$ \\ 
		\cline{4-7}
		&   & & DFP &$2m$ & 18 & $s=2$ \\ 
		\cline{2-7}
		& \multirow{3}{*}{2}  & \multirow{3}{*}{4} & UCP &$4m+20$& 45 & $s=4$ \\
		\cline{4-7}
		&   & & DFP &$4m$ & 25 & $s=3$\\ 
		\cline{4-7}
		&   & & HYP &$4m+2$  & 27 & $s=3$ \\ 
		\hline
		
	\end{tabular}
\end{table}
\clearpage

\addcontentsline{toc}{chapter}{References}
\renewcommand{\baselinestretch}{1}
\bibliographystyle{ieeetr}
\bibliography{BIBL}

\end{document}